\documentclass[11pt, reqno]{amsart}
\usepackage{amsmath,amsthm,amssymb}

\theoremstyle{plain}
\newtheorem{thm}{Theorem}[section]
\newtheorem{prop}[thm]{Proposition}
\newtheorem{lem}[thm]{Lemma}
\newtheorem{cor}[thm]{Corollary}

\theoremstyle{definition}

\newtheorem{rem}[thm]{Remark}
\newtheorem{defn}[thm]{Definition}
\newtheorem{eg}[thm]{Example}
\newtheorem{subtitle}[thm]{}
\newtheorem{ex}{Exercise}[section]
\numberwithin{equation}{section}

\def\a{\alpha}

\def\d{\delta}

\def\e{\epsilon}
\def\g{\gamma}

\def\K{\nabla}
\def\l{\lambda}

\def\n{\vert\/}
\def\o{\theta}

\def\W{\Omega}

\def\cv{{\mathcal{V}}}

\def\li{\langle}
\def\ri{\rangle}
\def\n{|\/ }
\def\tr{{\rm tr}}
\def\bs{\bigskip}
\def\ms{\medskip}
\def\ss{\smallskip}

\def\ti{\tilde}
\def\p{\partial}

\def\I{{\rm I\/}}

\def\diag{{\rm diag}}
\def\ad{{\rm ad}}

\def\C{\mathbb{C}}

\def\R{\mathbb{R} }
\def\Z{\mathbb{Z}}

\newcommand{\beq}{\begin{equation}}
\newcommand{\eeq}{\end{equation}}
\newcommand{\beg}{\begin{eg}}
\newcommand{\eeg}{\end{eg}}
\newcommand{\bthm}{\begin{thm}}
\newcommand{\ethm}{\end{thm}}
\newcommand{\bprop}{\begin{prop}}
\newcommand{\eprop}{\end{prop}}
\newcommand{\bcor}{\begin{cor}}
\newcommand{\ecor}{\end{cor}}
\newcommand{\blem}{\begin{lem}}
\newcommand{\elem}{\end{lem}}
\newcommand{\bca}{\begin{cases}}
\newcommand{\eca}{\end{cases}}
\newcommand{\brem}{\begin{rem}}
\newcommand{\erem}{\end{rem}}
\newcommand{\bpm}{\begin{pmatrix}}
\newcommand{\epm}{\end{pmatrix}}
\newcommand{\bbm}{\begin{bmatrix}}
\newcommand{\ebm}{\end{bmatrix}}
\newcommand{\bvm}{\begin{vmatrix}}
\newcommand{\evm}{\end{vmatrix}}
\newcommand{\bdefn}{\begin{defn}}
\newcommand{\edefn}{\end{defn}}
\newcommand{\bsub}{\begin{subtitle}}
\newcommand{\esub}{\end{subtitle}}
\newcommand{\bex}{\begin{ex}}
\newcommand{\eex}{\end{ex}}
\newcommand{\ben}{\begin{enumerate}}
\newcommand{\een}{\end{enumerate}}

\def\ti{\tilde}

\def\I{{\rm I\/}}

\def\diag{{\rm diag}}
\def\ad{{\rm ad}}

\def\rd{{\rm d\,}}
\def\id{{\rm id\/}}

\def\R{\mathbb{R} }
\def\C{\mathbb{C}}

\def\Z{\mathbb{Z}}

\def\calA{{\mathcal A}}

\def\calC{{\mathcal C}}
\def\calD{{\mathcal D}}
\def\calG{{\mathcal G}}
\def\calI{{\mathcal I}}
\def\calJ{{\mathcal J}}
\def\calK{{\mathcal K}}
\def\calL{{\mathcal L}}

\def\calO{{\mathcal O}}
\def\calP{{\mathcal P}}

\def\calU{{\mathcal U}}
\def\calV{{\mathcal V}}
\def\calX{{\mathcal X}}

\def\li{\langle}
\def\ri{\rangle}

\def\Res{{\rm Res\/}}

\begin{document}

\title[Tau functions and Virasoro Actions]
{Tau functions and Virasoro actions for soliton hierarchies}
\author{Chuu-Lian Terng$^\dag$}\thanks{$^\dag$Research supported
in  part by NSF Grant DMS-1109342}
\address{Department of Mathematics\\
University of California at Irvine, Irvine, CA 92697-3875.  Email: cterng@math.uci.edu}
\author{Karen Uhlenbeck$^*$}\thanks{$^*$Research supported in part by the Sid Richardson
Regents' Chair Funds, University of Texas system}
\address{The University of Texas at Austin\\ Department of Mathematics, RLM 8.100\\ Austin, TX 78712. Email:uhlen@math.utexas.edu}

\hskip 3in 

\begin{abstract} 
There is a general method for constructing a soliton hierarchy from a splitting $L_\pm$ of a loop group as a positive and a negative sub-groups together with a  commuting linearly independent sequence in the positive Lie algebra $\calL_+$.  Many known soliton hierarchies can be constructed this way. The formal inverse scattering associates to each $f$ in the negative subgroup $L_-$ a solution $u_f$ of the hierarchy.   When there is a $2$ co-cycle of the Lie algebra that vanishes on both sub-algebras,  Wilson constructed a tau function $\tau_f$ for each element $f\in L_-$.   In this paper, we give integral formulas for variations of $\ln\tau_f$ and second partials of $\ln\tau_f$, discuss whether we can recover solutions $u_f$ from $\tau_f$, and give a general construction of actions of the positive half of the Virasoro algebra on tau functions. We write down formulas relating tau functions and formal inverse scattering solutions and  the Virasoro vector fields for the $GL(n,\C)$-hierarchy. 
\end{abstract}

\maketitle

\section{Introduction}

This is the second in a series of papers attempting to give a uniform geometric structure in which many integrable systems can be placed. The two topics we address are tau functions and Virasoro actions, both of which are treated in a highly algebraic manner in the literature. Tau functions have a long and varied algebraic background, and both subjects appear as basic ingredients of several quantum cohomologies. The model theory is that of the KdV hierarchy, which generates the quantum cohomology of a point, and has been treated by many physicists and mathematicians, including Witten \cite{Wit91} and Konsevich \cite{Kon92}. 

Our approach is to define both the tau function and the Virasoro action on the space of solutions with the formal inverse scattering data in a simple conceptual way. Although the special solution of KdV for the quantum cohomology of a point does not have the right type of inverse scattering data, this approach gives all the correct formulae and allows us to define tau functions and Virasoro actions very generally. It also simplifies a number of computations in the literature.

To have suitable inverse scattering data, the hierarchies of interest must come from a (formal) Lie group $L$, and a splitting of the Lie algebra $\calL=\calL_+ + \calL_-$ with $L_+\cap L_-=\{e\}$.  The flows are generated by a sequence $J_j$ of commuting elements in $\calL_+$. The flows of the hierarchy are evolution partial differential equations on $C^\infty(\R, Y)$, where
$Y=[J_1,\calL_-]_+$ is the projection of $[J_1, \calL_-]$ to $\calL_+$ with respect to $\calL=\calL_+\oplus \calL_-$.

 Let 
\beq\label{kq1}
V(t)=\exp\left(\sum_{j=1}^N t_j J_j\right)
\eeq
denote the {\it vacuum frame\/}. 
 Given $f\in L_-$,  we can factor (cf. \cite{TerUhl11})
$$V(t)f^{-1}= M(t)^{-1} E(t)$$ with $M(t)\in L_-$ and $E(t)\in L_+$ for $t$ in an open subset of the origin in $\R^N$. Further
 $$u_f= (MJ_1M^{-1})_+ - J_1$$ 
 is a solution of the flows generated by $J_j$ for $1\leq j\leq N$. We call $u_f$ the {\it formal inverse scattering solution\/} given by scattering data $f$,  $E(t)$ the {\it frame\/} of $u_f$, and $M(t)$ the {\it reduced frame\/} for $u_f$.  We notice that the first flow equation for formal inverse scattering solutions is the translation, hence we can identify $x$ as $t_1$. 

The tau function is due to Wilson \cite{Wil91}. To construct it, one needs  a central extension {\it compatible with the splitting\/}, i.e., the $2$-cocyle for the central extension vanishes on $\calL_+$  and $\calL_-$. Then a function $\mu$ can be defined on the big cell $(L_+L_-)\cap (L_-L_+)$ using the difference in order of factorization.  This tau function is a local complex valued function $\tau_f$ of $t=(t_1, \ldots, t_N)$ and is defined to be $\mu(V(t) f^{-1})$ for scattering data $f\in L_-$. The construction of Wilson's tau function is simple and conceptual, but it is not easy to compute.  One result of this paper is to give integral formulas for derivatives and variations of $\ln \tau_f$ in terms of reduced frames of $u_f$.  Our formula for the first derivatives of $\ln\tau_f$ is similar to those which appear in many places including the work of Aratyn and van der Leur \cite{AraLeu03}.  
 
  The second derivatives of $\ln\tau_f$ are functions of the solution $u_f$. In ``good examples", solution $u_f$ can be recovered from the second partial derivatives of $\ln\tau_f$. For example, $u_f=-(\ln\tau_f)_{t_1t_1}$ for the KdV hierarchy.  
  Note that for many classical hierarchies including NLS, we cannot recover formal inverse scattering solutions $u_f$ for $\tau_f$. 
   
 We notice that a group acting on $L_-$ often produces a local action on the space of formal inverse scattering solutions and on the space of tau functions. Since we are interested in the Lie algebra of vector fields for these local group actions, these will automatically exist even though the group actions are not defined globally.  We would like to have nice formulae for them.  The formulas for the derivatives for tau functions and Virasoro actions appear in the literature without explanation. The example of $n\times n$ KdV is the most important, and is not given in this form elsewhere in the literature. 

We start by reviewing the definition of the Virasoro algebra. 
The Virasoro algebra $\cv$ is the Lie algebra spanned by $\{\xi_\ell\n \ell\in \Z\}$ with the bracket relations 
$$[\xi_j, \xi_k]= (k-j) \xi_{j+k}, \quad \forall\,\, j, k\in \Z.$$ The positive half Virasoro algebra $\cv_+$ is the sub-algebra of $\cv$ spanned by $\{\xi_j\n j\geq -1\}$.  An action of $\calV_+$ on a manifold is given by a sequence of tangent vector fields $X_j$ on the manifold satisfying the following bracket condition:
$$[X_j, X_k]=(k-j) X_{j+k}, \quad j, k\geq -1.$$ 
For example, $\calV_+$ acts on $S^1$ by 
$X_j(\l)= \l^{j+1}\frac{\p}{\p\l}$ for $ j\geq -1$.

Let $L(G)$ denote the group of smooth loops in a complex simple Lie group, and $L_\pm$ a splitting of $L(G)$ such that $L_+=L_+(G)$ is the subgroup of $f\in L(G)$ that is the boundary value of a holomorphic map from $|\l |<1$ to $G$.  Here is a simple recipe to construct Virasoro actions on the negative group $L_-$: Given a local group homomorphism $C$ from $S^1$ to $G$, 
\beq\label{km}
X_\ell(f)= -(\l^\ell(\l f_\l + fC'(1)) f^{-1})_-f, \quad \ell\geq -1,
\eeq
defines an action of $\calV_+$ on $L_-$, 
where $C'(1)=(\rd C/\rd \l)|_{\l=1}$.  
We can also give the formula for the induced $\calV_+$-action on $\ln\tau_f$. A ``good'' $\calV_+$-action on $\ln\tau_f$ should be given by partial differential operators of $\ln\tau_f$. To achieve this, we need to choose suitable homomorphism $C$ and then carry out long and complicated computations.

The $GL(n,\C)$ hierarchy is the hierarchy constructed form the standard splitting of the algebra of loops  in $GL(n,\C)$ generated by $\{a^i \l^j\n 1\leq i\leq n, j\geq 1\}$, where $a=\diag(c_1, \ldots, c_n)$ with distinct $c_i$'s. The flow equations are evolution equations on $C^\infty(\R, Y)$, where $Y=\{(\xi_{ij})\in gl(n,\C)\n \xi_{ii}=0, 1\leq i\leq n\}$. Let $A$ denote the subgroup of diagonal matrices in $GL(n,\C)$. For $k\in A$, we prove that $u_{kfk^{-1}}= k u_f k^{-1}$ but the second partials of $\ln\tau_f$ and $ \ln \tau_{kfk^{-1}}$ are equal. So we only expect to recover $A$ invariants of $u_f$ from $\ln\tau_f$. In fact, we show that 
$$(\ln\tau_f)_{t_{i,1} t_{j,1}}= (c_i-c_j)^2u_{ij}u_{ji}$$ for $1\leq i\not=j\leq n$, where $t_{i,1}= \sum_{i=1}^n c_i^j s_{j,1}$ and $s_{j,1}$ is the flow variable for the flow generated by $a^j\l$.   

We use $C(\l)=\I_n$ for the $GL(n,\C)$-hierarchy and prove that the Virasoro vector fields on $\calX=\ln\tau_f$ are given by partial differential operators of $\ln\tau_f$ (see Theorem \ref{do}).

The computations for the Virasoro actions on the tau functions and the proofs that we can recover formal inverse scattering solution $u_f$ from $\tau_f$ for the $n\times n$ KdV hierarchy are too complicated to include in this paper, and we leave these for the third paper \cite{TerUhl13b} in this series.
We prove in \cite{TerUhl13b} that we can recover $u_f$ from the second partials of $\ln\tau_f$ for the $n\times n$ KdV hierarchy.  We use $C(\l)= \diag(1,\l, \ldots, \l^{\frac{n-1}{n}})$ for the $n\times n$ KdV hierarchy and get the following Virasoro vector fields for 
$\calX= \ln\tau_f$, which agrees with that in the physics literature (\cite{vMo94}):
\begin{align*}
&\d_{-1}\calX =\frac{1}{n} \sum_{k>n}  k t_k \calX_{t_{k-n}} +\frac{1}{2n} \sum_{k=1}^{n-1} k(n-k) t_k t_{n-k},\\
&\d_0\calX = \frac{1}{n} \sum_{k\geq 1} k t_k \calX_{t_k},\\
&\d_\ell\calX = \frac{1}{n} \sum_{k\geq 1} kt_k \calX_{t_{n\ell+k}} + \frac{1}{2n}\sum_{k=1}^{n\ell-1} \left(\calX_{t_k}\calX_{t_{n\ell-k}} + \calX_{t_k t_{n\ell-k}}\right) \\
&\qquad \quad  +(\frac{1}{2n} -\frac{1}{2}) c_\ell(f), \quad \ell\geq 1,
\end{align*}
where $c_\ell(f)= \li \l^\ell ((\l f_\l  + f\Xi_0) f^{-1})^2\ri_0$, $\Xi_0= \frac{1}{n}\diag(0, 1, \ldots, (n-1))$, and $\calX_{t_{nk}}=0$. 

The outline of the paper is as follows: In section \ref{cd} we set up notation and review the construction and examples of integrable systems from splittings. In section \ref{eg}, we review the definition of Wilson's $\mu$ function dependent on a central extension compatible with the splitting and give integral formulas for variations of the $\mu$. The tau function is written in terms of the $\mu$-function. We write down formulas for second partials of $\ln\tau_f$ in terms of reduced frames,  give explicit relations between $\tau_f$ and $u_f$ in the the $2\times 2$ AKNS hierarchy and explain why we can not always recover $u_f$ from $\tau_f$ in section \ref{ki}.  We prove that \eqref{km} defines a $\calV_+$-action on $L_-$ and write down the induced $\calV_+$-action on $\ln\tau_f$ in terms of the reduced frame of $u_f$ in section \ref{kj}. Although we can not recover $u_f$ from second partials of $\ln\tau_f$ algebraically for the coupled $n$-component NLS hierarchy, we show that $u_f$ can be solved from an order $n$ linear system of ordinary differential equations in section \ref{ko}. In the last section, we compute the relation between $u_f$ and $\ln\tau_f$ and give explicit formulas of Virasoro vector fields for the $GL(n,\C)$-hierarchy.

\bs
\section{Soliton hierarchies constructed from Lie algebra splittings}\label{cd}

Here we set up notations and review the method of constructing soliton hierarchies from splittings of Lie algebras. For more details, see \cite{Ad79}, \cite{ReySem80}, \cite{DriSok84}, \cite{Sat84}, \cite{Wil91}, \cite{TerUhl00a}, \cite{Ter08}. 

\ss
Let $L_\pm$ be a splitting of the loop group $L=L(G)$, and  $\{ J_i\n i\geq 1\}$ a vacuum sequence (i.e., a linearly independent commuting sequence) in $\calL_+$. 
This data is enough, in good cases, to construct a hierarchy of commuting soliton flows on $C^\infty(\R, Y)$ of smooth maps from $\R$ to $Y$, where
\beq\label{dx}
Y=[J, \calL_-]_+.
\eeq

To get a good soliton theory, we must set things up so that the following two statements turn out to be true:

\ben
\item Given a smooth $u:\R\to Y$, there is a unique $Q(u)$ from $\R$ to $L$ satisfying 
\beq\label{ci}
\bca [\p_x-(J_1+u), \, Q(u)]=0,\\
Q(u)\, \, {\rm is \, \, conjugate \,\, to\,\,} J_1,\\
Q(0)=J_1,\eca
\eeq
Moreover, $Q(u)$ depends only on $u$ and $x$-derivatives of $u$.
\item For each $j\geq 1$, there is an analytic function $\phi_j$ and integers $s(j)$ such that
 \beq\label{jr}
 J_j=\phi_j(J_1)\l^{s(j)}.
 \eeq
 \een
 The second condition is not necessary, but it simplifies the computations in this paper.  For example, the $D_n$-KdV hierarchy constructed in \cite{DriSok84} does not satisfy \eqref{jr}. 
 
 \ms
 The flow in the hierarchy generated by $J_j$ is 
\beq\label{mb}
u_{t_j}= [\p_x-(J+u), (\phi_j(Q(u))\l^{s(j)})_+].
\eeq

\ms

\beg  \label{bq} {\bf The standard splitting of $L(G)$}\hfil\par

Let $L(G)$ denote the group of smooth maps $f$ from $S^1$ to a complex semi-simple Lie group $G$, $L_+(G)$ the subgroup of $f\in L(G)$ that can be extended holomorphically to $|\l|<1$, and $L_-(G)$ the subgroup of $f\in L(G)$ that can be extended holomorphically to $\infty\geq |\l|>1$ with $f(\infty)=\I$. Then $L_\pm(G)$  is called the {\it standard splitting\/} of $L(G)$. Let $\calG$ denote the Lie algebra of $G$, and $\calL(\calG)$ the Lie algebra of $L(G)$. The corresponding Lie subalgebras are
\begin{align*}
\calL_+(\calG)&=\{A\in \calL(\calG)\n A(\l)=\sum_{j\geq 0} A_j \l^j, \, A_j\in \calG\},\\
\calL_-(\calG)&=\{A\in \calL(\calG)\n A(\l)=\sum_{j<0}A_j \l^j, \, A_j\in \calG\}.
\end{align*}
\eeg

\ss
 The following hierarchies are the simplest known hierarchies given by splittings of $\calL(\calG)$ and subalgebras of $\calL(\calG)$ defined by finite order automorphisms of $\calG$.  

\beg\label{bi} {\bf The $G$-hierarchy}

Let $a_1, \ldots, a_n$ be a basis of a Cartan subalgebra $\calA$ of $\calG$ such that $a_1$ is {\it semi-simple\/}, i.e., the centralizer $\calG_{a_1}$ of $a_1$ in $\calG$ is $\calA$. 
The hierarchy constructed from the standard splitting $L_\pm(G)$ of $L(G)$ and the vacuum sequence 
$$\{J_{i,j}= a_i\l^j\n 1\leq i\leq n, j\geq 1\}$$ is the {\it $G$-hierarchy\/} on $C^\infty(\R, Y)$, where
$$Y=[J_{1,1}, \calL_-(\calG)]_+=\calG\cap \calA^\perp.$$
Here and henceforth we use $V^\perp$ to denote the orthogonal complement of a linear subspace $V$ of $\calG$ with respect to the Killing form on $\calG$. 
A direct computation implies that we can solve $Q(u)$ from \eqref{ci} (cf. \cite{Sat84}, \cite{TerUhl00a}) and $Q(u)= a_1\l + u +\sum_{j<0} Q_j(u) \l^j$. Moreover, the flow generated by $J_{1,1}$ is $u_{t_1}= u_x$ and the flow generated by $J_{i,1}= a_i \l$ for $2\leq i\leq n$ is the {\it $n$-wave equation\/}
$$u_{t_{i,1}}= \ad(a_i)\ad(a_1)^{-1}(u_x) -[u, \ad(a_i)\ad(a_1)^{-1}(u)].$$

For $\calG= sl(2,\C)$, we let $a=\diag(i, -i)$. 
The $SL(2,\C)$-hierarchy constructed from $L_\pm(SL(2,\C))$ and the vacuum sequence $\{a\l^j\n j\geq 1\}$ is the {\it $2\times 2$ AKNS hierarchy\/}.  The flows are evolution equations for
$$u=\bpm 0& q\\ r &0\epm.$$  Write 
$Q(u)=a\l + \sum_{i\leq 0} Q_i(u) \l^i$ and solve $Q(u)$ from \eqref{ci}. We obtain 
$$Q_{-1}(u)= \frac{i}{2}\bpm qr& -q_x\\ r_x & -qr\epm, \quad Q_{-2}(u)= \frac{1}{4}\bpm -qr_x+ q_xr & -q_{xx} + 2q^2 r\\ -r_{xx} + 2qr^2 & qr_x - q_x r \epm.$$
So the flow equations \eqref{mb} generated by $J_j=a\l^j$ with $j=1, 2,3$ are
\begin{align*}
&q_{t_1}= q_x, \quad r_{t_1}= r_x,\\
&q_{t_2} = -\frac{i}{2}(q_{xx} -2q^2r), \quad r_{t_2}= \frac{i}{2} (r_{xx} - 2q r^2),\\
& q_{t_3}= -\frac{1}{4}(q_{xxx} - 6qq_x r), \quad r_{t_3}= \frac{1}{4} (r_{xxx} - 6 qrr_x).
\end{align*}
\eeg

\beg \label{bk} {\bf The $U$-hierarchy}\par

Let $\tau$ be a group involution of $G$ such that the differential of $\tau$ at the identity $e$ (still denoted by $\tau$) is conjugate linear (i.e., $\tau(c \xi)= \bar c\tau(\xi)$ for all $c\in \C$ and $\xi\in \calG$). The fixed point set $U$ of $\tau$ is  a {\it real form of $G$\/}. Let $L^\tau(G)$ denote the subgroup $g\in L(G)$ satisfying the {\it $U$-reality condition\/}
$$\tau(g(\bar\l))=g(\l),$$
and $L^\tau_\pm(G)= L^\tau(G)\cap L_\pm(G)$.  The Lie algebra is 
$$\calL^\tau(\calG)=\{\sum_i \xi_i \l^i\,\n\, \xi_i \in \calU\}.$$
Let $\{a_1, a_2, \ldots, a_n\}$ be a basis of a maximal abelian sub-algebra $\calA_\R=\calA \cap \calU$ such that the centralizer $\calU_{a_1}$ of $a_1$ in $\calU$ is $\calA_\R$.  The hierarchy constructed from the splitting $L^\tau_\pm(U)$ and the vacuum sequence $\{J_{i,j}= a_i\l^j\n 1\leq i\leq n, j\geq 1\}$ is called the {\it $U$-hierarchy\/}.  

   For example, let $G=SL(2,\C)$, $\tau(g)=(\bar g^t)^{-1}$, and $a_1=\diag(i, -i)$. Then we have $U=SU(2)$ and the flow generated by $a\l^2$ is the {\it non-linear Schr\"odinger equation\/} (NLS)
$$r_t= \frac{i}{2}(r_{xx} + 2|r|^2 r).$$ 
So  the $SU(2)$-hierarchy is the {\it NLS hierarchy\/}. 

  \eeg
  
  \beg\label{bkc} {\bf The $(G,\sigma)$-hierarchy}\par

Let $\sigma$ be a group involution of $G$, whose differential at $e$ (still denoted by $\sigma$) is complex linear.  Let $L^\sigma(G)$ denote the subgroup of $g\in L(G)$ such that $g$ satisfies the {\it $\sigma$-reality condition\/}:
$$g(\l)=\sigma(g(-\l)),$$
and $L^\sigma_\pm(G)= L_\pm(G)\cap L^\sigma(G)$.  Decompose $\calG= \calG_0\oplus \calG_1$ as $\pm 1$ eigenspaces of $\sigma$ on $\calG$. Then 
$\sum_i \xi_i \l^i \in \calL^\sigma(\calG)$ if an only if $\xi_i\in \calG_0$ for $i$ even and $\xi_i\in \calG_1$ for $i$ odd. 
 Let $\{a_1, \ldots, a_m\}$ be a basis of a maximal abelian sub-algebra $\calA_1$ in $\calG_1$ such that $\ad(a_1)$ maps $\calA_1^\perp$ isomorphically onto $(\calG_0)_a^\perp$. The hierarchy constructed from the splitting $L_\pm^\sigma(G)$ and the vacuum sequence $\{a_i\l^{2j-1}\n 1\leq i\leq m, j\geq 1\}$ is called the {\it $(G,\sigma)$-hierarchy\/}.  

For example, let $G=SL(2,\C)$, $c= \bpm 0&1\\ 1&0\epm$, and $\sigma$ the involution on $SL(2,\C)$ defined by $\sigma(g)= cgc^{-1}$. Let $\{a\l^{2j-1}\n j\geq 1\}$ be the vacuum sequence, where $a=\diag(1,-1)\in \calG_1$. Then the third flow in the $(SL(2,\C), \sigma)$-hierarchy is the {\it complex mKdV\/}:
$$q_t=\frac{1}{4}(q_{xxx} - 6q^2 q_x).$$

\eeg

\beg\label{bla} {\bf The $\frac{U}{K}$-hierarchy}\par 

Let $\tau$ and $\sigma$ be commuting group involutions of $G$ such that the induced Lie algebra homomorphisms $\tau$ is conjugate linear and $\sigma$ is complex linear.  Set 
$$L^{\tau, \sigma}(G)= L^\sigma(G)\cap L^\tau(G), \quad L_\pm^{\tau, \sigma}(G)= L_\pm(G)\cap L^{\tau, \sigma}(G).$$
Note that $f\in L^{\tau,\sigma}(G)$ if and only if $f$ satisfies the following $\frac{U}{K}$-hierarchy:
$$\tau(f(\bar\l))=f(\l), \quad \sigma(f(-\l))= f(\l).$$
Then $L_\pm^{\tau,\sigma}(G)$ is a splitting of $L^{\tau,\sigma}(G)$.  
Let $U$ denote the fixed point set of $\tau$ in $G$, and $K$ the fixed point set of $\sigma$ in $U$. Then $\frac{U}{K}$ is a symmetric space. Let  $\calU=\calK\oplus \calP$ be the Cartan decomposition for the symmetric space, i.e., $\calK$ and $\calP$ are $+1, -1$ eigen-spaces of $\sigma$ on $\calU$.  Then $\sum_i \xi_i\l^i\in \calL^{\tau,\sigma}(\calG)$ if and only if $\xi_i\in \calK$ for even $i$ and $\xi_i\in \calP$ for odd $i$. 
Let $\calA_0$ be a maximal abelian sub-algebra in $\calP$, and $\{a_1, \ldots, a_n\}$ a basis of $\calA_0$ such that $\ad(a_1)$ maps the orthogonal complement of $\calA_0^\perp$ in $\calP$ isomorphically onto $\calK_{\calA_0}^\perp\cap \calK$.  The hierarchy constructed from the splitting $L_\pm^{\tau, \sigma}(G)$ and the vacuum sequence $\{a_i\l^{2j-1}\n 1\leq i\leq n, j\geq 1\}$ is called the {\it $\frac{U}{K}$-hierarchy\/}. 

 For example, let $G=SL(2,\C)$, $\tau(g)= (\bar g^t)^{-1}$, $\sigma(g)=(\bar g^t)^{-1}$, and $a_1=\diag(i, -i)$. Then $U=SU(2)$, $K=SO(2)$, and flows in the $\frac{SU(2)}{SO(2)}$-hierarchy are evolution equations for maps 
$$u=\bpm 0&-r\\ r&0\epm, \quad r\in C^\infty(\R, \R).$$
Moreover, the flow generated by $a\l^3$ is 
 the mKdV equation
$$r_t=\frac{1}{4} (r_{xxx} + 6 r^2 r_x).$$ 
\eeg

Suppose $g\in L(G)$ and $g=g_+g_-$ with $g_\pm \in L_\pm(G)$.  It can be easily seen that if $g\in L^\tau(G)$ then $g_\pm \in L^\tau_\pm(G)$ and similar statements are true for $g$ in $L^\sigma(G)$  and $L^{\tau,\sigma}(G)$. 
Since the formal inverse scattering solution $u_f$ of the $G$-hierarchy is constructed from factorization of $V(t)f^{-1}$ for $f\in L_-(G)$, we obtain the following proposition.

\bprop\label{ie} Let $u_f$ be the formal inverse scattering solution defined by $f\in L_-(G)$. Then $u_f$ is a solution of
\ben
\item the $U$-hierarchy if $f\in L_-^\tau(G)$,
\item the $(G,\sigma)$-hierarchy if $f\in L_-^\sigma(G)$,
\item the $\frac{U}{K}$-hierarchy if $f\in L_-^{\tau,\sigma}(G)$.
\een 
\eprop

\beg\label{cz} {\it The KdV-hierarchy}

It can be shown that the condition $q=1$ is invariant under the flows generated by $a\l^{2j+1}$ for all $j\geq 0$ in the $2\times 2$ AKNS hierarchy. The restriction gives the KdV hierarchy  (cf. \cite{AKNS74}).  
Unlike the NLS and mKdV, we do not know the condition on $f\in L_-(SL(2,\C))$ such that the solution $u_f$ of the $2\times 2$ AKNS hierarchy satisfies the constraint $q=1$ for the KdV. 
However, 
the KdV hierarchy can also be constructed from an unusual splitting of $L(SL(2,\C))$ as  follows (cf. \cite{TerUhl11}):  
 Let $L^{kdv}$ denote the subgroup of $f\in L(SL(2,\C))$ satisfying
\beq\label{bz}
\phi(\l) f(\l) \phi(\l)^{-1}= \phi(-\l) f(-\l) \phi(-\l)^{-1},
\eeq
where 
$$ \phi(\l)=\bpm 1&0\\ \l&1\epm.$$
Then $L^{kdv}_\pm= L^{kdv}\cap L_\pm(SL(2,\C))$ is a splitting of $L^{kdv}$.  
Let $J= \bpm \l & 1\\ 0&-\l\epm$. The hierarchy constructed from the splitting $L^{kdv}_\pm$ and the vacuum sequence $\{J^{2j-1}\n j\geq 1\}$ is the KdV hierarchy.  
Moreover, a direct computation implies that the the solution $Q(u)$ of \eqref{ci} for $u=\bpm 0&0\\r&0\epm$ is of the form
$$Q(u)= J+ \sum_{i\leq 0} Q_i(u)\l^i= a\l + e_{12} +  u + \sum_{i<0} Q_i \l^i $$
with 
$$Q_{-1}= -\frac{1}{2} \bpm r&0\\ r_x &-r\epm, \quad Q_{-2}= \frac{1}{4} \bpm r_x & -2r\\ r_{xx}- 2r^2 & -r_x\epm.$$
So the flow generated by $J^3$ is the KdV equation
$$r_t= \frac{1}{4}(r_{xxx} - 6r r_x).$$
\eeg

A direct computation gives the following formula for the solution of \eqref{ci} for $u=u_f$:

\bprop\label{fv}
 Let $L_\pm$ be a splitting of $L$, $\{J_j\n j\geq 1\}$ a vacuum sequence, $V(t)$ the vacuum frame defined by \eqref{kq1}, and $f\in L_-$. Let $u_f$ denote the formal inverse scattering solution given by $f$, and $M$ the reduced frame of $u_f$, i.e., 
$V(t)f^{-1}= M(t)^{-1} E(t)$
with $M(t)\in L_-$ and $E(t)\in L_+$ and $u_f= (MJ_1M^{-1})_+-J_1$. Then 
$$[\p_x+J_1+ u_f, MJ_1M^{-1}]=0,$$ 
i.e., $Q(u_f)= MJ_1M^{-1}$ is the solution of \eqref{ci}  for $u=u_f$.
\eprop

\brem The map from $f\in L_-$ to the formal inverse scattering solution $u_f$ is not injective. In fact, 
if $h\in L_-$ commutes with $J_1$, then $u_{fh}= u_f$. To see this, we factor 
$V(t)(fh)^{-1}=\ti M^{-1}(t) \ti E(t)$ with $\ti M(t)\in L_-$ and $\ti E(t)\in L_+$. But $V(t)h^{-1} f^{-1}= h^{-1}V(t) f^{-1}= h^{-1} M^{-1}E$ and $Mh\in L_-$ imply that $\ti M= Mh$. Hence $MJ_1M^{-1}= \ti M J_1 \ti M^{-1}$ and $u_{fh}= u_f$. 
\erem

\bs
\section{Formulas for variations of Wilson's $\mu$ functions}\label{eg}

We review Wilson's construction of function $\mu$ for splittings and give integral formulas for the variation of $\mu$.  Wilson's tau function is defined in terms of the $\mu$-function. We find it surprising the function $\mu$ defined by Wilson \cite{Wil91} is not better known.  

Start with an infinite dimensional Lie group $L$ with a splitting  $L_\pm$ of $L$. The {\it big cell\/} $\calC$ of $L$ is
$$\calC=(L_+L_-)\cap (L_-L_+).$$  
Hence an element $g\in L$ lies in $\calC$ if and only if $g$ can be factored as $f_+f_-$ and also as $g_-g_+$ uniquely with $f_\pm, g_\pm\in L_\pm$. 
In general, $\calC$ is open and contains the identity of $L$, and in good cases, it will be dense. Since we are only interested in local formulas in this paper, the Local Factorization Theorem (cf. \cite{TerUhl11}) is enough, i.e., if $\g$ is a smooth map from an open neighborhood $\calO$ of the origin in $\R^k$ to $L$ and $\g(0)$ lies in the big cell $\calC$, then there is an open subset $\calO_0$ in $\R^k$ containing the origin such that $\g(t)\in \calC$ for all $t\in \calO_0$. 

First we recall the definition of a central extension of a Lie algebra. 
A $2$-cocycle of a Lie algebra $\calL$ is a skew-symmetric bilinear form $w$ on $\calL$ satisfying 
$$w([\xi_1, \xi_2], \xi_3)+ w([\xi_2, \xi_3], \xi_1) + w([\xi_3, \xi_1], \xi_2)=0$$
for all $\xi_i\in \calL$.  

Define $[\, ,]_1$ on 
$$\hat \calL= \calL+ \C c$$
by
$$[\xi + rc, \eta]_1= [\xi, \eta] + w(\xi, \eta) c$$
for $\xi, \eta\in \calL$ and $r\in \C$.  Then 
$$\C c\to \hat \calL \to \calL$$ is a central extension of $\calL$, where the projection $\pi: \hat \calL\to \calL$ is defined by $\pi(\xi+ rc)=\xi$. 

We choose the $2$-cocycle $w$ such that the left invariant form on $L$ defined by $w$ (still denoted by $w$) is an integral cohomology class, i.e., $w\in H^2(L, 2\pi i\Z)$. Let 
$$\C^*\to \hat L \to L$$
denote the central extension of the group $L$, i.e., a principal $\C^\ast$-bundle whose Chern class is $w$. Here $\C^*=\C\setminus \{0\}$. Next we review the  construction of $\hat L$  in Pressley and Segal \cite{PreSeg86}, which we need to use to construct natural lifts of $L_\pm$ to $\hat L$.   Let
 $$P\times \C^*:=\{(\g, z)\n g:[0,1]\to L \, \, {\rm such \, that \, } \g(0)=e, z\in \C^*\}$$
 denote the product space of smooth paths in $L$ starting at the identity  and $\C^*$. 
We define $(\g_1, z_1)\sim (\g_2, z_2)$ if $\g_1(1)=\g_2(1)$ and
$$z_2= z_1 \exp\left(\int_{\W(\g_1, \g_2)} w\right),$$
where $\W(\g_1, \g_2)$ is a surface bounded by $\g_1\ast I(\g_2)$. Here $\g_1\ast \g_2$ denote the composition of paths (the path $\g_1$ followed by $\g_2$) and $I(\g)(s)= \g(1-s)$ is the path $\g$ with reverse orientation. The condition that $w$ is an integral cohomology class implies that this is a well-defined equivalence relation.  The central extension 
$$\hat L=(P\times \C^*)/\sim$$ is the group of equivalence classes $[(\g, z)]$ and the multiplication is given by
$$[(\g_1, z_1)]\cdot [(\g_2, z_2)] = [(\g_1\ast(g_1 \g_2), z_1z_2)], \quad {\rm where\,} \, g_1=\g_1(1).$$
The projection $\pi: \hat L\to L$ defined by $\pi([(\g, z)])= \g(1)$ is a principal $\C^*$-bundle with first Chern class $w$.

\bdefn The central extension $\hat L$ of $L$ constructed from the $2$-cocycle $w$ is {\it compatible\/} with the splitting $L_\pm$ of $L$ if $w$ vanishes on both $\calL_+$ and $\calL_-$.  Or equivalently, $\calL_\pm$ are isotropic subspaces of $\calL$. \edefn

\bprop\label{kr}
There is a natural lift $S:L_+\cup L_-\to \hat L$ with $\pi\circ S=\id$.
\eprop

\begin{proof}
Define $S(g_\pm) = [(\g_\pm, 1)]$, where $\g_\pm$ is any path in $L_\pm$ joining $e$ to $g_\pm$.  Since $w$ vanishes on $L_\pm$, $S$ is independent of the choice of path $\g_\pm$ as long as it lies in $L_\pm$.
\end{proof}

  Note that the restriction of   $S$ to $L_\pm$ are group homomorphisms, i.e., 
$$S(g_1g_2)= S(g_1)S(g_2)$$ if both $g_1, g_2$ lie in $L_+$ or both lie in $L_-$. 

\bdefn \label{fa} The {\it $\mu$ function\/} is a complex function on the big cell $\calC= (L_+L_-)\cap (L_-L_+)$ defined as follows: Given $f\in \calC$, then there exist unique $f_\pm, g_\pm \in L_\pm$ such that 
$$f= f_+f_-^{-1}= g_-^{-1} g_+.$$
Let $S$ denote the natural lift of $L_+\cup L_-$ to $\hat L$ defined by Proposition \ref{kr}.
Since $S(f_+)S(f_-^{-1})$ and $S(g_-^{-1})S(g_+)$ lie in the same fiber of the principal $\C^*$-bundle $\hat L$ over $L$, they differ by a scalar in $\C^*$, which we call $\mu(f)$. In other words, $\mu(f)$ is defined by the following identity:
$$S(f_+)S(f_-^{-1})= \mu(f) S(g_-^{-1})S(g_+).$$
\edefn

To compute $\mu(f)$, we first factor $f=f_+f_-^{-1}= g_-^{-1}g_+$ with $f_\pm, g_\pm \in L_\pm$.  Let $\g_\pm$ be a path joining $e$ to $f_\pm$ in $L_\pm$, and $\ti \g_\pm$ a path in $L_\pm$ joining $e$ to $g_\pm$.  Then $\mu(f)$ is given by the integral 
\beq\label{ac} \mu(f)= \exp\left(\int_\W w\right),\eeq
where $\W$ is any surface bounded by the curves $\g_+\ast (f_+\g_-^{-1})$ and $\ti \g_-^{-1}\ast(g_-^{-1} \ti \g_+)$.  

   \bthm {\bf Variations of the $\mu$ function} \label{cn} \hfil
  
Let $w$ be an integral $2$-cocycle on $L$ compatible with the splitting $L_\pm$ of $L$, $B$ the unit disk in $\R^2$, and $h:(-r, r)\times B\to L$ a smooth map. Set $h_\e= h(\e, \cdot)$.  Then 
\beq\label{hd}
\frac{d}{d\e} \int_B h_\e^* w=\int_{\p B} w(h^{-1}\p_\e h, h^{-1}\p_sh) ds.
\eeq
\ethm

\begin{proof} Let $(x,y)$ denote coordinates of $B$.  Then 
$$\int_B h_\e^*w= \int_B w(h^{-1}h_x\, , h^{-1}h_y) dx dy.$$
So 
\beq\label{dt}
\frac{d}{d\e} \int_B h_\e^*w =\int_B w((h^{-1}h_x)_\e, h^{-1}h_y) + w(h^{-1}h_x, (h^{-1}h_y)_\e) dx dy.
\eeq
Use the condition that $w$ is a 2 co-cycle and 
\begin{align*}
(h^{-1}h_x)_\e &= (h^{-1}h_\e)_x + [h^{-1}h_x, h^{-1}h_\e],\\
(h^{-1}h_y)_\e &= (h^{-1}h_\e)_y + [h^{-1}h_y, h^{-1}h_\e],
\end{align*} 
to imply that the integrand of the right hand side of \eqref{dt} is equal to
\beq\label{du}
w((h^{-1}h_x)_\e, h^{-1}h_y) + w(h^{-1}h_x, (h^{-1}h_\e)_y) + w([h^{-1}h_x, h^{-1}h_y], h^{-1}h_\e).
\eeq
Let $\o$ denote the $1$-form on $B$ defined by
$$\o(X)= w(h^{-1}h_\e, X).$$
By Cartan's formula, we have $\rd \o (X,Y)= X\o(Y) - Y\o(X) - \o([X,Y])$, So $\rd \o(h^{-1}h_x,h^{-1}h_y)$ is equal to \eqref{du}. Then formula \eqref{hd} follows from Stokes' Theorem.
\end{proof}

\ss
\bsub {\bf A $2$ co-cycle on the loop algebra $\calL(\calG)$}

Let $G$ be a Lie group, $\calG$ the Lie algebra of $G$,  $L(G)$ the group of  smooth loops $f:S^1\to G$, and $\calL(\calG)=C^\infty(S^1, \calG)$ the Lie algebra of $L(G)$.
Assume that $(\, , )$ is an ad-invariant non-degenerate bilinear form on $\calG$.  For a classical Lie algebra $\calG$,  
$$(A, B )= \tr(AB),$$ is such a bilinear form.
Fix an integer $k$. Given $\xi, \eta\in \calL(\calG)$, define 
\beq \label{aa}\li \xi, \eta\ri_k  =\, {\rm the\, coefficient\, of\, } \l^{k} \, {\rm of\, } (\xi(\l),\eta(\l)) =\sum_j \tr(\xi_j \eta_{-j+k}),
\eeq
i.e., $(\xi(\l), \eta(\l))= \sum_k \li\xi, \eta\ri_k \l^k$. 
Then
\beq\label{ab} 
w(\xi, \eta)=\li \p_\l \xi, \eta\ri_{-1}= \sum_j j\tr(\xi_j \eta_{-j})
\eeq
is a $2$-cocycle on $\calL(\calG)$ (cf.  \cite{PreSeg86}). 
\esub

We will apply the following theorem to $W(t)= V(t) =\exp(\sum_{j=1}^n t_j J_j)$ in section \ref{ki}. 

\bthm\label{df}  Let $\li\, ,\ri_{-1}$ denote the bilinear form on $\calL(\calG))$ defined by \eqref{aa}, and let $L_\pm$ be a splitting of $L(G)$ such that $\li \calL_+, \calL_+\ri_{-1} =0$ and the co-cycle $w$ defined by \eqref{ab} vanishes on $\calL_\pm$, i.e., $w$ is compatible with the splitting. Let $\mu$ be the function defined on the big cell $\calC= (L_+L_-)\cap (L_-L_+)$ as in Definition \ref{fa}, 
$W(t)$ a path in  $L_+$,  and $f_-\in L_-$. If
$$W(t)f_-^{-1}= M(t)^{-1}E(t), \quad {\rm with \/}\,\, M(t)\in L_- \, {\rm and\/}\, E(t)\in L_+,$$   then 
\beq\label{eb}
\frac{d}{d\, t} \ln \mu(W(t)f_-^{-1}) = \li (\p_t W)W^{-1}, M^{-1}\p_\l M\ri_{-1}.
\eeq
\ethm

\begin{proof}
Let $p(s)$ be a path in $L_-$ that joining $e$ to $f_-$, $\ti W(t, s)= W(st)$, $\ti M(t, s)$ a path  in $L_-$ joining $e$ to $M(t)$,  and $\ti E(t, s)$ a path in $L_+$ joining $e$ to $E(t)$ for each $t$.  Let $\li\, , \ri$ denote $\li\, , \ri_{-1}$ in the calculation below. 

We use \eqref{hd} to compute $\frac{d}{dt}\ln (\mu(W(t)f^{-1}))$. 
The boundary has four pieces: $\ti W(t, s)$, $W(t)p(s)^{-1}$, $M(t)^{-1} \ti E(t, s)$, and $\ti M(t, s)^{-1}$.  Since $\ti W(t,s)\in L_+$ and $\ti M(t, s)\in L_-$ and $i_\pm^*w=0$,  the boundary integral vanishes on the first and last pieces. Then
\eqref{hd} gives 
\beq\label{hda}
\frac{d}{d\, t} \ln \mu(W(t)f_-^{-1}) =\int_0^1 w(\g^{-1}\g_t, \g^{-1}\g_s)\, ds- \int_0^1 w(h^{-1}h_t, h^{-1}h_s)\, ds.
\eeq
Here $\g(t,s)= W(t) p^{-1}(s)$ and $h(t,s)= M(t)^{-1} \ti E(t, s)$. 
Note that 
\begin{align*}
&\g^{-1}\g_t= pW^{-1}W_tp^{-1}, \quad \g^{-1}\g_s= -p_sp^{-1},\\
& h^{-1}h_t= -\ti E^{-1} M_t M^{-1} \ti E + \ti E^{-1} \ti E_t, \quad h^{-1}h_s= \ti E^{-1} \ti E_s,
\end{align*}
and $i_\pm^*w =0$.  Hence the right hand side of \eqref{hda} is equal to 
\begin{align*}
&\int_0^1 w(pW^{-1}W_t p^{-1}, -p_sp^{-1}) ds - \int_0^1 w(\ti E^{-1} M_tM^{-1}\ti E, \ti E^{-1} \ti E_s) ds\\ 
&=\int_0^1 \li pW^{-1}W_tp^{-1}, (p_s p^{-1})_\l\ri ds -\int_0^1 \li \ti E^{-1} M_tM^{-1} \ti E, (\ti E^{-1}\ti E_s)_\l\ri ds\\
&= \int_0^1 \li W^{-1}W_t, p^{-1}(p_sp^{-1})_\l p\ri ds -\int_0^1 \li M_tM^{-1}, \ti E(\ti E^{-1}\ti E_s)_\l \ti E^{-1}  \ri ds.
\end{align*}
A direct computation gives 
$$p^{-1}(p_sp^{-1})_\l p= (p^{-1}p_\l)_s, \quad \ti E(\ti E^{-1}\ti E_s)_\l \ti E^{-1}= (\ti E_\l \ti E^{-1})_s.$$
We use $\li\, , \ri$ to denote $\li\, ,\ri_{-1}$ in the rest of the proof.
The above integral is equal to
\begin{align*} 
&\li W^{-1}W_t, p^{-1}p_\l\ri\, \big|_{s=0}^{s=1} - \li M_tM^{-1}, \ti E_\l \ti E^{-1}\ri\, \big|_{s=0}^{s=1}\\
&= \li W^{-1}W_t , f^{-1} f_\l\ri - \li M_t M^{-1}, E_\l E^{-1}\ri.
\end{align*}
Use $M=EfW^{-1}$ and $\li \calL_+, \calL_+\ri =0$ to compute the second term:
\begin{align*}
& \li M_tM^{-1}, E_\l E^{-1}\ri = \li E_t E^{-1}- EfW^{-1}W_t W^{-1} M^{-1}, E_\l E^{-1}\ri\\
& = -\li EfW^{-1}W_t W^{-1}M^{-1}, E_\l E^{-1}\ri 
=-\li MW_tW^{-1} M, E_\l E^{-1}\ri\\
&= -\li MW_tW^{-1}M, M_\l M^{-1} + MW_\l W^{-1} M^{-1} - MW f^{-1}f_\l W^{-1}M^{-1}\ri\\
&= -\li W_t W^{-1}, M^{-1}M_\l\ri + \li W^{-1}W_t, f^{-1} f_\l\ri.
\end{align*}
This proves the theorem. 
\end{proof}

\bs
\section{The partial derivatives of tau functions}\label{ki}

We use the variation formula for $\mu$ given in Theorem \ref{df} to compute the partial derivatives of $\ln\tau_f$ in flow variables. We then prove that the second partials of $\ln\tau_f$ are polynomial in $u_f$ and its $t_1$ derivatives.  We also show that there is a finite dimensional symmetry on the negative group $L_-(G)$ that leaves $\tau_f$ invariant but acts on $u_f$ non-trivially for the $G$-hierarchy.  This explains why we can not recover $u_f$ from $\tau_f$ for the NLS hierarchy.  Such a symmetry does not exist for KdV, which has a different $J_1$.  

First we recall the definition of $\tau_f$ given by Wilson in \cite{Wil91}.

\bdefn\label{kf} Assume that $L_\pm$ is a splitting of $L$ compatible with the $2$-cocyle that defines a central extension $\C^*\to \hat L\to L$, and that $\calJ=\{J_j\n j\geq 1\}$ is a vacuum sequence in $\calL_+$.
For $f\in L_-$, the {\it tau function $\tau_f$ associated to $f$\/} is a function of $t=(t_1, \ldots, t_N)$ defined by 
$$\tau_f(t)= \mu(V(t) f^{-1}),$$ 
where $V(t)=\exp(\sum_{j=1}^N J_j t_j)$ is the vacuum frame and $\mu$ is the $\mu$-function defined in Definition \ref{fa}. 
\edefn

\brem \label{br} Let $L_\pm$ be a splitting of the loop group $L(G)$, and $V(t)$ the vacuum frame.  Given  $f\in L_-$, we have $V(0)f^{-1}= f^{-1}$ is in $L_-$, which is contained in the big cell. So it follows from the Local Factorization Theorem (Theorem 1.2 of \cite{TerUhl11}) that given any $f\in L_-$, there exists an open neighborhood $\calO_0$ of the origin in $\R^N$ such that $V(t)f^{-1}$ lies in the big-cell $\calC$ for all $t\in \calO_0$.  Hence  $\tau_f(t)$ is defined for  $t\in \calO_0$. 
\erem

Next we use \eqref{eb} to calculate the derivatives of $\ln \tau_f$.

\bthm\label{cs} Let $L_\pm\subset L$ be a splitting, $\calJ=\{J_j\n j\geq 1\}$ a vacuum sequence, $w$ a $2$-cocycle on $\calL$ compatible with the splitting, and $V(t)=\exp(\sum_{j=1}^N t_j J_j)$ the vacuum frame.  Let $f\in L_-$, and  
$$V(t)f^{-1}= M^{-1}(t) E(t)$$ with $M(t)\in L_-$ and $E(t)\in L_+$.   Then 
\ben
\item $(\ln \tau_f)_{t_j}= \li J_j, M^{-1}\p_\l M\ri_{-1} = \li MJ_jM^{-1}, (\p_\l M)M^{-1}\ri_{-1}$,
\item $(\ln \tau_f)_{t_1t_j}= \li MJ_j M^{-1}, \p_\l J_1\ri_{-1}$,
\een
where $\li\, , \ri_{-1}$ is the bilinear form defined by \eqref{aa}.
\ethm

\begin{proof}
Since $V^{-1}\p_{t_j}V= J_j$, (1) follows from \eqref{eb}.  

Recall that we have $(\p_{t_1}M)M^{-1}= -(MJ_1M^{-1})_-$ and $(MJ_1M^{-1})_+= J_1+u$.  In the proof below, we use $\li\, ,\ri$ to denote $\li\, , \ri_{-1}$. By (1), we have
\begin{align*}
&\p_{t_1}\p_{t_j}\ln \tau= \li J_j, \p_{t_1} (M^{-1} \p_\l M)\ri\\
&= \li M J^jM^{-1}, M\p_{t_1} (M^{-1} \p_\l M) M^{-1}\ri \\
&= \li MJ_j M^{-1}, \p_\l((\p_{t_1}M) M^{-1})\ri\\
& =\li MJ^jM^{-1}, \p_\l(-(MJ_1M^{-1})_-\ri,
\end{align*}
and since $\li MJ_jM^{-1}, \p_\l(MJ_1M^{-1})\ri=0$, the last equality is equal to 
\begin{align*}
&= \li MJ_jM^{-1}, \p_\l(MJM^{-1}-(MJ_1M^{-1})_-)\ri \\
&= \li MJ_jM^{-1}, \frac{\p}{\p \l} (MJ_1M^{-1})_+\ri\\
&= \li MJ_jM^{-1}, \p_\l(J_1+u)\ri = \li MJ_jM^{-1}, \p_\l J_1\ri.  
\end{align*}
\end{proof} 

Note that formulae for $(\ln\tau_f)_{t_j}$ appears in the literature in many places giving as definition of tau functions, in particular in the work of Aratyn and van der Leur \cite{AraLeu03}. They proved that $\Theta= \sum_{j=1}^\infty \li MJ_j M^{-1}, (\p_\l M)M^{-1}\ri_{-1} \rd t_j$ is a closed $1$-form, hence it is $d$ of a function, which they call $\ln\tau_f$.  So their $\ln\tau_f$ can differ by a constant with Wilson's $\ln\tau_f$. 

We use similar computations as for Theorem \ref{cs} to get the following.

\bthm\label{csk} With the same assumption as in Theorem \ref{cs}. Then 
\beq\label{if}
(\ln \tau)_{t_j t_k}= \li MJ_jM^{-1}, \p_\l(MJ_kM^{-1})_+ \ri_{-1}.\eeq
\ethm

Suppose  $\li \calL_+, \calL_+\ri_{-1} = \li \calL_-, \calL_-\ri_{-1} =0$. Below we compute directly to see that the right hand side of \eqref{if} is symmetric in $j$ and $k$, what it should be true from Theorem \ref{csk}.   We use $\li\, , \ri$ to denote $\li\, ,\ri_{-1}$. A direct computation implies that
\begin{align*}
&\li MJ_j M^{-1}, \p_\l (MJ_k M^{-1})_+\ri = \li (MJ_j M^{-1})_-, \p_\l (MJ_k M^{-1})_+\ri\\
&= - \li \p_\l(MJ_j M^{-1})_-, (MJ_kM^{-1})_+\ri = -\li \p_\l(MJ_j M^{-1})_-, MJ_k M^{-1}\ri \\
&= -\li \p_\l (MJ_j M^{-1})-\p_\l (MJ_jM^{-1})_+, MJ_k M^{-1}\ri\\
&= -\li \p_\l (MJ_j M^{-1}), MJ_kM^{-1}\ri + \li \p_\l (MJ_jM^{-1})_+, MJ_k M^{-1}\ri\\
&= -\li [M_\l M^{-1}, MJ_jM^{-1}], MJ_kM^{-1}\ri + \li \p_\l (MJ_jM^{-1})_+, MJ_k M^{-1}\ri\\
&=-\li M_\l M^{-1}, [MJ_jM^{-1}, MJ_kM^{-1}]\ri + \li \p_\l (MJ_jM^{-1})_+, MJ_k M^{-1}\ri\\
&= \li \p_\l (MJ_jM^{-1})_+, MJ_k M^{-1}\ri.
\end{align*}
In the last equality we used the facts that $\li\, , \ri$ is ad-invariant and $[J_j, J_k]=0$. 

\brem
Given $h\in (L_-)_{J_1}$ and $f\in L_-$, we have seen that $u_f= u_{fh}$. What is the relation between $\tau_f$ and $\tau_{fh}$? By Theorem \ref{cs}, we have  
$$(\ln\tau_{fh})_{t_j}= (\ln\tau_f)_{t_j} + \li J_j, h_\l h^{-1}\ri_{-1}.$$
In other words, $(\ln \tau_f)_{t_j}$ and $(\ln\tau_{fh})_{t_j}$ differ by some constant $c_j$ independent of $t$ for each $j\geq 1$.
This implies that $(\ln\tau_f)_{t_it_j}= (\ln\tau_{fh})_{t_i t_j}$. 
\erem

 By Proposition \ref{fv}, the solution $Q(u_f)$ of \eqref{ci} is $MJ_1M^{-1}$, where $M$ is   the reduced frame of $u_f$.  Recall that we assume $J_j= \phi_j(J_1)\l^{s(j)}$ for some analytic function $\phi_j$ and some non-negative integer $s(j)$. Therefore we have $MJ_jM^{-1}= \phi_j(MJ_1M^{-1}))\l^{s(j)}$.  It follows from Theorem \ref{csk} and the fact that $Q(u_f)$ depends only on $u_f$ and its $t_1$ or $x$ derivatives that we have the following corollary.
 
\bcor
$\ln(\tau_f)_{t_jt_k}$ is a function of $u_f$ and its $t_1$-derivatives. 
\ecor

\ms

\beg\label{cj}  We give explicit formulas of $(\ln\tau_f)_{t_1t_j}$ in terms of $u_f$ for the $SL(2,\C)$-hierarchy and its various restrictions. 

Let $a=\diag(1,-1)$, and $f\in L_-(SL(2,\C))$. Write $u_f=  \bpm 0& q\\ r &0\epm$. By Proposition \ref{fv}, we have $Q(u_f)= MJM^{-1}$. Write $Q(u_f)$ in power series in $\l$:
$$Q(u_f) = MJM^{-1}=a\l + \sum_{i\leq 0} Q_i\l^i.$$ 
 By Theorem \ref{cs} (2),  we have
$$(\ln\tau_f)_{t_1t_j}=\tr(aQ_{-j}).$$
Use the formulas of $Q_i$ given in Example \ref{bi} to get
\beq\label{bh}
\bca
(\ln\tau_f)_{t_1t_1}= \tr(aQ_{-1})= -qr, \\
 (\ln\tau_f)_{t_1t_2}=\tr(aQ_{-2}) = \frac{1}{2}(q_{t_1} r- r_{t_1}q).\eca
 \eeq
 Hence $\ln\tau_f$ does not determines $u_f$.  
 But a simple computation implies that
 \beq\label{kg}
 \bca q_{t_1}= -(\frac{y_2}{y_1}+\frac{ (y_1)_{t_1}}{2y_1})\, q,\\
 r_{t_1}= (\frac{y_2}{y_1}-\frac{(y_1)_{t_1}}{2y_1})\, r,\eca 
 \eeq
 where $y_1= (\ln\tau_f)_{t_1t_1}$ and $y_2= (\ln \tau_f)_{t_1 t_2}$. 
 This shows that $u_f$ is related to $\ln\tau_f$ by a system of first order linear ordinary differential equations, whose coefficients are rational functions of $(\ln\tau_f)_{t_1t_1}$ and $(\ln\tau_f)_{t_1 t_2}$.    

The first restriction gives the NLS hierarchy.  
 If $f\in L_-(SL(2,\C))$ satisfies the $SU(2)$-reality condition, $\overline{f(\bar\l))}^t f(\l)=\I$, then $u_f=\bpm 0& q\\ -\bar q &0\epm$ and $q$ is a solution of the NLS hierarchy. 
By \eqref{bh}, we have 
$$(\ln\tau_f)_{t_1t_1}= |q|^2, \quad (\ln\tau_f)_{t_1t_2}= -\frac{1}{2}(q_{t_1}\bar q- q \bar q_{t_1}).$$
Write $q= \rho e^{i\o}$ in polar coordinates. Then 
$$(\ln\tau_f)_{t_1t_1}= \rho^2, \quad (\ln\tau_f)_{t_1 t_2}= - \o_{t_1} \rho^2.$$
Hence $\ln\tau_f$ determines $u_f$ up to a constant in $S^1$.  

We also obtain the mKdV-hierarchy with different reality conditions. 
If $f\in L_-(SL(2,\C))$ satisfies the $\frac{SU(2)}{SO(2)}$-reality condition,
$$\overline{f(\bar\l))}^t f(\l)=\I, \quad f(-\l)= (f(\l)^t)^{-1},$$
then $u_f= \bpm 0& q\\ -q &0\epm$ is a solution of the mKdV hierarchy. 
 By Theorem \ref{cs}(2),  
$$(\ln\tau_f)_{t_1t_1}= -q^2.$$ 
So $\ln\tau_f$ determines $u_f$ up to a sign.
\eeg

\beg {\it Tau functions for the KdV hierarchy}

We have seen in Example \ref{cz} that KdV can be obtained in two ways.  We will calculate tau functions in both ways. 
 The flows of $2\times 2$ AKNS hierarchy generated by $a\l^{2j-1}$  leave the condition $q=1$ invariant and the resulting odd flows give the KdV hierarchy.  Use $Q_{-1}$ given in Example \ref{bi} with $q=1$ to get  $Q_{-1}=\frac{i}{2}\bpm r & 0\\ r_x & -r\epm$.   By Theorem \ref{cs}(2) we have 
 $$(\ln\tau_f)_{t_1t_1}= -r.$$ 
 
 The KdV hierarchy also can be constructed from the splitting of $\calL^{kdv}$ and vacuum sequence $\{J^{2j-1}\n j\geq 1\}$ as given in Example \ref{cz}.   Apply Theorem \ref{cs} to see that again $(\ln\tau_f)_{t_1t_1}= -r$. 
 \eeg
 
 The vector AKNS is a natural generalization of the $2\times 2$ AKNS hierarchy. We will prove that the formal inverse scattering solution $u_f$ can be solved from a system of linear ordinary differential equations from $\ln\tau_f$. Since the proof is complicated, we will do this in section \ref{ko}.  
 
\ms
 
The next theorem gives a natural finite dimension group action on $L_-(G)$ for the standard splitting. The induced action on second partials of $\ln\tau_f$ is trivial but the induced action on the formal inverse scattering solution for the $G$-hierarchy is non-trivial. This explains why we can not recover $u_f$ from $\tau_f$.

\bthm\label{fb}
Let $L=L(G)$,  $L_\pm=L_\pm(G)$ the standard splitting of $L$,   and $\{J_j\n j\geq 1\}$ a vacuum sequence with $J_1=a\l$ for some regular $a\in \calG$.  Let  $G_a$ denote the subgroup of $k\in G$ that commutes with $a$. Then for  $k\in G_a$ and $f\in L_-$ we have  
\ben
\item[(i)]  $kfk^{-1}\in L_-$,
\item[(ii)] $u_{kfk^{-1}}= ku_f k^{-1}$,
\item[(iii)]  $(\ln\tau_f)_{t_1t_j}= (\ln\tau_{kfk^{-1}})_{t_1t_j}$ for all $j\geq 1$,
\een
where $u_f$ and $u_{kfk^{-1}}$ are solutions constructed from scattering data $f$ and $kfk^{-1}$ respectively.
\ethm

\begin{proof}
It is clear that $kL_\pm k^{-1}= L_\pm$.  Let $V(t)$ denote the vacuum frame. Set $\ti f= kfk^{-1}$. Then $\ti f\in L_-$.  Factor $V(t)f^{-1}= M(t)^{-1}E(t)$ and $V(t)\ti f^{-1}= \ti M(t)^{-1}\ti E(t)$  with $M(t), \ti M(t)\in L_-$ and $E(t), \ti E(t)\in L_+$.  Since $kJ_1=J_1k$, $kV(t)=V(t)k$.  By assumption, we have $kL_\pm k^{-1}\subset L_\pm$.  So we have $\ti M= kMk^{-1}$ and $\ti E= kEk^{-1}$. But $E_{t_1}E^{-1}= a\l+ u_f$ and $\ti E_{t_1}\ti E^{-1}= a\l + u_{\ti f}$.  Hence $u_{\ti f}= ku_f k^{-1}$. 

We use Theorem \ref{cs}(2) and the fact that $k$ commute with $J_j$ and $a$ to compute
\begin{align*}
(\ln\tau_{kfk^{-1}})_{t_1t_j} &= \li \ti M J_j \ti M^{-1}, a\ri =\li kMk^{-1}J_j kM^{-1}k^{-1}, a\ri \\
&=\li kMJ_jM^{-1}k^{-1}, a\ri = \li MJ_jM^{-1}, k^{-1}a k\ri\\
&= \li MJ_jM^{-1}, a\ri = (\ln \tau_f)_{t_1t_j}.
\end{align*}
This proves Statement (3).  
\end{proof}

As a consequence, we only expect to recover $G_a$-invariants of $u_f$ from second partials of $\ln\tau_f$.  

\ms
The proof of Theorem \ref{fb} also implies the next result. 

\bthm\label{fba}
Let $G, \tau, \sigma$, $U$, $K$, $G_0$ and $a=a_1$ be as in section \ref{cd} for the $G$, $U$, $(G,\sigma)$, and $\frac{U}{K}$ hierarchies.   
\ben
\item If $k\in U_a$ and $f\in L^\tau_-(G)$, then $kfk^{-1}\in L^\tau_-(G)$.  
\item If $k\in K_a$ and $f\in L^{\tau, \sigma}_-(G)$, then $kfk^{-1}\in L^{\tau, \sigma}_-(G)$.
\item If $k\in (G_0)_a$ and $f\in L_-^\sigma(G)$, then $kfk^{-1}\in L_-^\sigma(G)$.
\een
Moreover, (ii) and (iii) of Theorem \ref{fb} hold for the $U$-, $(G,\sigma)$-, and the $\frac{U}{K}$-hierarchies. 
\ethm

\beg \hfil \par

(1) For the $2\times 2$ AKNS hierarchy, we have $a=\diag(i, -i)$. So
$$K=SL(2,\C)_a=\{\diag(c, c^{-1})\n c\in \C, c\not=0\}.$$   Given $f\in L_-$ and $k\in K$, it follows from Proposition \ref{fb} that we have
 $(\ln\tau_f)_{t_1 t_j}= (\ln \tau_{kfk^{-1}})_{t_1 t_j}$ and $u_{kfk^{-1}}= k u_f k^{-1}$.  Write $u_f=\bpm 0& q\\ r&0\epm$ and $k=\diag(c, c^{-1})$ for some non-zero $c\in \C$. Then $ku_f k^{-1}= \bpm 0 & c^2 q\\ c^{-2} r &0\epm$. So we can only recover $K$-invariants of $u_f$ from $\ln\tau_f$ (note that the right hand side of \eqref{bh} is invariant under the action of $K$).

(2) For the NLS hierarchy, we have  $$SU(2)_a=\{\diag(\a, \a^{-1})\n \a\in \C, |\a |=1\}\simeq S^1$$ and the action of $SU(2)_a$ on $q$ is $\diag(\a, \a^{-1})\ast q= \a^2 q$ for $\a\in S^1$.  By Proposition \ref{fb}, $(\ln\tau_f)_{t_1t_j}= (\ln \tau_{kfk^{-1}})_{t_1t_j}$ for $j=1,2$.  Note that both $|q|$ and $\o_{t_1}$ are invariant under the action of $S^1$. 

(3) For the mKdV hierarchy, we have $SO(2)_a= \Z_2$ and $q^2$ is invariant under the $Z_2$-action. 
\eeg

\bs

\section{Virasoro actions}\label{kj}

In this section on Virasoro actions, we derive the infinitesimal formulas for vector fields using local group actions. It is, of course true, that once we know the formulas, we can compute directly the formulas for Lie brackets and invariance of the flows. However, these are long, abstract, and somewhat unmotivated calculations.  Hence this section approaches the Virasoro vector fields as infinitesimal flows derived from local group actions, the formulas for group actions are more direct, better motivated, and easier to understand. 

Let $L_\pm$ be a splitting of the group $L(G)$ of smooth loops in the complex Lie group $G$ such that $L_+=L_+(G)$, the subgroup of $g\in L(G)$ that is the boundary value of a holomorphic map from $|\l |\leq1$ to $G$ (this is the same $L_+$ as in the standard splitting).   In this section, we first associate to each local homomorphism from $S^1$ to $G$ an action of the positive half Virasoro algebra $\calV_+=\{\xi_j\n j\geq -1\}$ on $L_-$.  Then we compute the induced action on reduced frames and the tau functions for the hierarchies given by the splitting $L_\pm$ and a vacuum sequence $\{J_j\n j\geq 1\}$. 

We use the following simple fact a number of times.  

 \bprop\label{cia}
 Let $L_\pm$ be a splitting of $L$,  $g$ a smooth curve on $L$ defined on an open interval $(-r, r)$ for some $r>0$, $g(0)\in L_-$, and $g(s)= g_+(s) g_-(s)$ with $g_\pm\in L_\pm$. Let $'$ denote the $s$ derivative.  Then 
 \beq\label{cl}
g_-'(0) g(0)^{-1}= (g'(0) g(0)^{-1})_-.
 \eeq
 \eprop
 
 \begin{proof}  
  Since $g(0)\in L_-$, we have $g_+(0)=\I$.  Lemma follows from a simple computation. 
  \end{proof}
 
 Let $D_+(S^1)$ denote the group of diffeomorphisms of $S^1$ that is the boundary value of a holomorphic map from $|\l | \leq 1$ to itself. Recall that $\calV_+$ acts on $S^1$ and  $\calV_+$ can be embedded as a sub-algebra of the Lie algebra $\calD_+(S^1)$ by $\{X_j= \l^{j+1} \frac{\p}{\p_\l}\n j\geq -1\}$, where $\l=e^{i\o}$. 
 
Let $C$ be a local group homomorphism from $S^1$ to $G$, i.e., $C$ is defined on $\calI_\e=\{e^{i\o}\,\n\, |\o |< \e\}$ for some $\e>0$ such that $C(\l_1\l_2)= C(\l_1)C(\l_2)$ whenever $\l_1, \l_2, \l_1\l_2\in \calI_\e$. It is helpful to know that the basic example is $C(\l)\equiv\I_n$ the trivial homomorphism. The Theorem below associate to each $C$ a natural $\calV_+$-action on $L_-$.

\bthm\label{eu}
Let $L_\pm$ be a splitting of $L=L(G)$ with  $L_+=L_+(G)$, and  $C$ a local group homomorphism from $S^1$ to $G$. Given $k\in D_+(S^1)$ and $f\in L_-$, define 
\begin{subequations}
\begin{gather}
(k\diamond f)(\l)= f(k^{-1}(\l)) C\left(\frac{k^{-1}(\l)}{\l}\right), \label{as}\\
(k\sharp f)= {\rm the \,} L_- \, {\rm component\,\, of \,\,} k\diamond f, \label{at}
\end{gather}
\end{subequations}
i.e.,  $k\sharp f= g_-$ is obtained from factoring $k\diamond f= g_+g_-$ with $g_\pm \in L_\pm$.  
Then $\sharp$ defines an action of $D_+(S^1)$ on $L_-$.
Moreover, the infinitesimal vector field $Z_j$ on $L_-$ corresponding to $X_j=\l^{j+1} \frac{\p}{\p\l}$ for the $\sharp$-action is 
\beq\label{bb}
Z_j(f) f^{-1}= -(\l^{j+1} f_\l f^{-1}+ \l^j fC'(1) f^{-1})_-, \quad {\rm for\,\, } j\geq -1,
\eeq
where $C'(1)= \frac{\rd C}{\rd \l}|_{\l=1}$.  
In particular, \eqref{bb} defines an action of $\calV_+$ on $L_-$.
\ethm

\begin{proof} First note that if $f_+\in L_+$ and $k\in D_+(S^1)$ then $f_+\circ k\in L_+$. 

Let $k_1, k_2\in D_+(S^1)$, $(k_1\circ k_2)\sharp f=\ti g_-$, $k_2\sharp f= g_-$, and $k_1\sharp g_-= h_-$.  So we have
\begin{subequations}
\begin{gather}
f(k_2^{-1}(k_1^{-1}(\l)) C\left(\frac{k_2^{-1}(k_1^{-1}(\l))}{\l}\right) = \ti g_+(\l) \ti g_-(\l),\label{ba1}\\
f(k_2^{-1}(\l)) C\left(\frac{k_2^{-1}(\l)}{\l}\right)= g_+(\l) g_-(\l), \label{ba2}
\end{gather}
\end{subequations}
where $g_\pm$, $h_\pm$, and $\ti g_\pm$ are in $L_\pm$.  

We use \eqref{ba1} and \eqref{ba2} and the assumption that $C$ is a group homomorphism to compute $k_1\sharp g_-$:
\begin{align*}
& g_-(k_1^{-1}(\l)) C\left(\frac{k_1^{-1}(\l)}{\l}\right)\\
&= g_+^{-1}(k_1^{-1}(\l)) f(k_2^{-1}(k_1^{-1}(\l))C\left(\frac{k_2^{-1}(k_1^{-1}(\l)}{k_1^{-1}(\l)}\right)C\left(\frac{k_1^{-1}(\l)}{\l}\right)\\
&=  g_+^{-1}(k_1^{-1}(\l)) f(k_2^{-1}(k_1^{-1}(\l))C\left(\frac{k_2^{-1}(k_1^{-1}(\l)}{\l}\right)\\
&= g_+^{-1}(k_1^{-1}(\l)) \ti g_+(\l) \ti g_-(\l).
\end{align*}
But $g_+^{-1}\circ k_1^{-1}\in L_+$. So $\ti g_-$ is the $L_-$ factor of $k_1\diamond g_-$, which implies that $\ti g_-= k_1\sharp g_-= k_1\sharp (k_2\sharp f)$. 

Use Proposition \ref{cia} and a straight forward computation to get the formula for the infinitesimal vector fields of the action $\sharp$ on $L_-$ corresponding to $L_j$.
\end{proof}
   
\beg\label{ay} Choose $C(\l)=\I$ in Theorem \ref{eu}. Then 
\beq\label{aj}
\zeta_j(f)= -(\l^{j+1} \frac{\p f}{\p \l} f^{-1})_- f, \quad j\geq -1.
\eeq
defines an action of $\calV_+$ on $L_-$.
\eeg

\beg\label{az} For $G=L(GL(n,\C))$, choose  
$$C(\l)= \diag(1, \l^{1/n}, \ldots, \l^{(n-1)/n})$$ (a local homomorphism). 
By Theorem \ref{eu}, we  get the following $\calV_+$-action on $L_-$,
\beq\label{ap}
\d_j(f)= -(\l^{j+1} \frac{\p f}{\p \l}f^{-1} + \l^j f\Xi_0f^{-1})_- f, \quad j\geq -1,
\eeq
where 
\beq\label{ap1}
\Xi_0:=C'(1)= \frac{1}{n} \diag(0, 1, \ldots, n-1).
\eeq

\eeg

Next we give the formula for the induced variation on the reduced frame when we vary $f\in L_-$.

\bprop \label{dj}
Let $L_\pm$ be a splitting of $L$, $\{J_j\n j\geq 1\}$ a vacuum sequence in $\calL_+$, and $V(t)=\exp(\sum_{i=1}^N t_i J^i)$.  Let $\d f$ denote a variation of $f$ in $L_-$, and $\d M$ and $\d E$ the corresponding variations of $M$ and $E$ respectively computed from $V(t)f^{-1}= M^{-1}(t) E(t)$.  Then 
\beq\label{bf}
(\d M) M^{-1}= (E(\d f) f^{-1} E^{-1})_-.
\eeq
\eprop

\begin{proof}
Use $E=MVf^{-1}$ and a direct computation to get
$$(\d E) E^{-1}= (\d M) M^{-1} -E(\d f)f^{-1} E^{-1}.$$
Since $\d E$ is a variation on $L_+$,  $(\d E) E^{-1}\in \calL_+$ and the Proposition follows. 
\end{proof}

The following Theorem gives the variation of $\tau_f$ when we vary $f$.

\bthm\label{by} Let $L_\pm$ be a splitting of $L(G)$ with $L_+=L_+(G)$, $\{J_j\n j\geq 1\}$ a vacuum sequence in $\calL_+$, and $V(t)=\exp(\sum_{j=1}^N t_j J_j)$ the vacuum frame. Let $f(\e)$ be a curve in $L_-$. Factor $V(t)f(\e)^{-1}= M^{-1}(t, \e) E(t, \e)$ with $E(t, \e)\in L_+$ and $M(t,\e)\in L_-$.  Then 
\beq\label{ah} 
\frac{\p}{\p \e} \ln\tau_{f(\e)}= -\li M_\e M^{-1}, E_\l E^{-1}\ri_{-1},  
\eeq 
and $M_\e M^{-1}= (Ef_\e f^{-1} E^{-1})_-$, where $M_\e=\p M/\p \e$ and $E_\l= \p E/ \p\l$. In other words, the induced variation on $\ln\tau$ when we vary $f$ is given by
\beq\label{ah1}
\d \ln\tau_f= -\li \d M M^{-1}, E_\l E^{-1}\ri_{-1}.
\eeq
\ethm

\begin{proof}
Fix $\e$ and $t$, let $p(s,\e)$ be a path in $L_-$ that joins $e$ to $f(\e)$, $\ti V(t, s)= V(st)$, $\ti M(t, s, \e)$ a path  in $L_-$ joining $e$ to $M(t, \e)$,  and $\ti E(t, s,\e)$ a path in $L_+$ joining $e$ to $E(t,\e)$.  
We use \eqref{hd} and proceed similarly as in the proof of Theorem \ref{df} to compute $\frac{\p}{\p\e} \ln\mu(V(t)f(\e))$. We get
\beq\label{hdb}
\frac{\p}{\p\e} \ln \mu(V(t)f(\e)^{-1}) =\int_0^1 w(\g^{-1}\g_\e, \g^{-1}\g_s)\, ds- \int_0^1 w(h^{-1}h_\e, h^{-1}h_s)\, ds,
\eeq
where $\g(t,s,\e)= V(t) p^{-1}(s,\e)$ and $h(t,s,\e)= M(t,\e)^{-1} \ti E(t, s,\e)$. 
A direct computation gives 
\begin{align*}
&\g^{-1}\g_\e= -p_\e p^{-1}, \quad \g^{-1}\g_s= -p_sp^{-1},\\
& h^{-1}h_\e= -\ti E^{-1} M_\e M^{-1} \ti E + \ti E^{-1} \ti E_\e, \quad h^{-1}h_s= \ti E^{-1} \ti E_s.
\end{align*}
Since the 2-cocycle vanishes on $\calL_\pm$ and $\ti E(\ti E^{-1}\ti E_s)_\l \ti E^{-1}= (\ti E_\l \ti E^{-1})_s$, we have
\begin{align*}
&\frac{\p}{\p\e} \ln \mu(V(t)f(\e)^{-1})= -\int_0^1\li \ti E^{-1} M_\e M^{-1}\ti E, (\ti E^{-1} \ti E_s)_\l\ri\\
&\quad =-\int_0^1\li M_\e M^{-1} \ti E(\ti E^{-1}\ti E_s)_\l \ti E\ri =-\int_0^1\li M_\e M^{-1}, (\ti E_\l \ti E^{-1})_s\ri\\
& \quad = -\li M_\e M^{-1}, \ti E_\l \ti E^{-1})|_{s=0}^{s=1}= -\li M_\e M^{-1}, E_\l E^{-1}\ri,
\end{align*}
Notice the boundary term at $s=0$ is zero because $\ti E(t, s,\e)= e$ when $s=0$. This completes the proof of \eqref{ah}.  

Since $M(t,\e)=E(t, \e) f(\e) V(t)^{-1}$, a direct computation implies that
$$M_\e M^{-1}= E_\e E^{-1} + Ef_\e f^{-1} E^{-1}.$$
But the left hand side lies in $\calL_-$ and $E_\e E^{-1}\in \calL_+$, hence we get the formula for $M_\e M^{-1}$. 
\end{proof}

From \eqref{bf} and \eqref{ah1}, we get the following theorem.

\bthm \label{id} Let $L_\pm$ and $C$ be as in Theorem \ref{eu}, and $\{Z_\ell\n \ell\geq -1\}$ the $\calV_+$-action on $L_-$ defined by \eqref{bb}. Then the induced $\calV_+$-action on reduced frames and tau-functions are
\begin{align}
\d_\ell (M) M^{-1} &= -(\l^\ell E(\l f_\l f^{-1}+ fC'(1) f^{-1})E^{-1})_-,\label{kp1}\\
\d_\ell(\ln\tau) &= \li \l^\ell E(\l f_\l f^{-1} + fC'(1) f^{-1}) E^{-1}, \l E_\l E^{-1}\ri_0.\label{kp2}
\end{align} 
\ethm

Ideally, we would like the induced action on $\ln\tau_f$ given by differential operators in $\ln\tau_f$.  But formula \eqref{kp2} is written in terms of the frame $E=MVf^{-1}$. Although we have formulas to express derivatives of $\ln\tau_f$ in terms of $M$, it is not clear whether these vector fields \eqref{kp2} on $\ln\tau_f$ are differential operators in $\ln\tau_f$. We will prove that this is the case for the $GL(n,\C)$-hierarchy in section \ref{pc} when we choose $C=\I_n$. We will prove in the third paper \cite{TerUhl13b} of our series that this is the case for the $n\times n$ KdV hierarchy if we use $C$ as in Example \ref{az}. We also prove in  \cite{TerUhl13b}  that there is a diffeomorphism between the phase spaces of the $n\times n$ KdV hierarchy  and  the GD$_n$-hierarchy (the Gelfand-Dikki hierarchy on the space of order $n$ linear differential operators on the line) such that the flows correspond and that our $\calV_+$-action on $\ln\tau_f$ is the known $\calV_+$-action  under this diffeomorphism. 

\bs 
\section{Tau functions for vector AKNS  hierarchy}\label{ko}

We have seen in Example \ref{cj} that although we can not recover $u_f$ from second partials of $\ln\tau_f$ for the $2\times 2$ AKNS-hierarchy, we can solve $u_f$ from a linear system of ordinary differential equations with coefficients being rational functions of  second partials of $\ln\tau_f$. We show that this is true for a natural generalization of the $2\times 2$ AKNS hierarchy that having the vector NLS hierarchy as a restriction. 

For this system, we choose the standard splitting $\calL_\pm$ of $\calL=\calL(gl(n+1,\C))$,
 choose $\{a\l^j\n j\geq 1\}$ to be the vacuum sequence, where 
 \beq\label{js}
 a=\bpm i\I_n&0\\ 0 & -i\epm.
 \eeq
Then the phase space is $C^\infty(\R, Y)$, where 
$$Y= [J_1, gl(n+1,\C)]_+= \left\{ u=\bpm 0& q\\ r&0\epm\,\big|\, r, q^t\in \C^{n\times 1}\right\}.$$
We call this hierarchy the {\it  vector AKNS\/}. It generalizes the $2\times 2$ AKNS, where $r, q\in \C$.  The further restriction $$r=-\bar q^t$$ gives the hierarchy (cf. \cite{ForKul83}), whose second flow is the {\it coupled $n$-component NLS\/} or the {\it vector NLS\/}
$$q_t= \frac{i}{2}(q_{xx} + 2||q||^2 q).$$
 This hierarchy has the important property of being gauge equivalent to the hierarchy belonging to Schr\"odinger flow from $\R\times \R\to \C P^n$ (cf. \cite{TerUhl06}). 
The further restriction to the case when $u=\bpm 0& -r^t\\ r&0\epm\in gl(n+1,\R)$ gives hierarchy containing the vector mKdV,
$$r_t= \frac{1}{4}(r_{xxx}+ 6||r||^2 r_x).$$ The vector mKdV is the principal curvature flow of the {\it geometric Airy curve flow on $\R^{n+1}$}, 
$$\g_t= \K_{e_1}^\perp H(\g),$$
where $H(\g)$ is the mean curvature vector of the curve $\g(\cdot, t)$ in $\R^{n+1}$, $e_1$ is the unit tangent to $\g(\cdot, t)$, and $\K^\perp$ is the induced normal connection (\cite{Ter14a}).

We have seen from Theorem \ref{cs} that 
$$\{y=(y_1, \ldots, y_{2n})\n y_j= (\ln\tau_f)_{t_1t_j}\n 1\leq j\leq 2n\}$$
are $2n$ complex differential polynomials of $(r, q^t)\in C^{2n}$ in $t_1$ variable.  Our goal is to recover $u_f$ from these $y_i$'s by solving a system of linear ordinary differential equations. 

 Write
$$Q(u_f)=MJ_1M^{-1}= a\l + u_f + \sum_{j<0} Q_j \l^j.$$ Since $J_j= J_1\l^{j-1}$, we get $MJ_jM^{-1}= Q(u_f) \l^{j-1}$.  By Theorem \ref{cs}, we have
\beq\label{ia}
(\ln\tau_f)_{t_1 t_j}= \li MJ_j M^{-1}, a\ri_{-1}= (Q_{-j}, a).
\eeq
To obtain a more precise relation between 
$\{y_j=(\ln\tau_f)_{t_1t_j}\n 1\leq j\leq 2n\}$ 
 and $u_f$, we need to compute $Q(u_f)$ from \eqref{ci}. In particular, $Q(u_f)$ satisfies 
\beq\label{ck}
\bca [\p_x-(a\l+u_f), \, Q(u_f)]=0,\\ Q(u_f)^2= -\l^2 \I_{n+1}. \eca
\eeq 
First we set up some notation: Let
\begin{align*}
\calG_0:&= gl(n+1,\C)_a=\left\{\bpm \xi &0\\ 0& c\epm\,\bigg|\, \xi\in gl(n,\C), c\in \C\right\},\\
 \calG_1:&= (gl(n+1,\C)_a)^\perp=\left\{\bpm 0 &\eta_2\\ \eta_1 &0\epm\,\bigg|\, \eta_1, \eta_2^t\in \C^{n\times 1}\right\}.
 \end{align*} 
 Then $gl(n+1,\C)= \calG_0\oplus \calG_1$, 
 $[\calG_0, \calG_0]\subset \calG_0$, $[\calG_0, \calG_1]\subset \calG_1$, and $[\calG_1, \calG_1]\subset \calG_0$.
We write $u^{(k)}= \p_x^k u$.
Define the {\it weight $\nu$\/} for monomials as follows:
$$\nu(u^{(i_1)} u^{(i_2)} \cdots u^{(i_k)})= \sum_{m=1}^k (i_m+1).$$

\bprop\label{cm}
  Let $Q(u)=a\l + u +\sum_{j<0} Q_j (u)\l^j$
  be the solution of \eqref{ck}, and  write $Q_{-j}(u)= P_{-j}(u) + T_{-j}(u)$ with $P_{-j}(u)\in\calG_1 $ and $T_{-j}(u)\in \calG_0$. Then  
\begin{align}
&P_{-j}(u)= (-a/2)^j u^{(j)} + \phi_j,\label{co1}\\
&T_{-j}(u)=(a/2)^j\left(\sum_{i=0}^{j-1} (-1)^i u^{(i)} u^{(j-1-i)} \right) + \psi_{j-1}, \label{co2}
\end{align}
where $\phi_j$ is a  weight $(j+1)$ polynomial in $u^{(i)}$ and $ au^{(i)}$ for $ i\leq j-2$ and $\psi_j$ is a weight $(j+1)$ polynomial in $u^{(i)}$ and $au^{(i)}$ for $i\leq j-1$.  
\eprop

\begin{proof}
A simple calculation shows that for all $y\in \calG_1$ and our choice of $a$, we have
\beq\label{cp}
ay=-ya, \quad [a,y]= 2ay, \quad a^{-1}=-a.
\eeq
To simplify the notation, we write $Q=Q(u)$ and $Q_j= Q_j(u)$. 
Equate coefficients of $\l^{-j}$ in \eqref{ck} to get
$$(Q_{-j})_x-[u, Q_{-j}]= [a, Q_{-(j+1)}].$$
Equate the $\calG_0$ and $\calG_1$ components of the above equation to get 
 \begin{subequations}
 \begin{gather}
P_{-(j+1)} = -\frac{1}{2} a((P_{-j})_x + [T_{-j}, u], \label{aj1}\\
-2aT_{-(j+1)} = \sum_{i=0}^{j} P_{-i} P_{-(j-i)} + T_{-i} T_{-(j-i)}.\label{aj2}
\end{gather}
\end{subequations}
Then the lemma follows from induction on $j$. 
\end{proof}
 
\ms

For example, 
$T_0=\phi_0=\phi_1= \psi_0=0$, and
\begin{align*}
&Q_1=a, \quad Q_0=u, \quad Q_{-1}= \frac{a}{2}(-u_x+ u^2),\\
&Q_{-2}= - \frac{1}{4} u_{xx} + \frac{1}{2} u^3 - \frac{1}{4} (uu_x - u_x u).
\end{align*}

Since
\beq\label{cx} J_j= J_1\l^{j-1}= a\l^j,
\eeq
the flows \eqref{mb} generated by $J_j$ in the vector AKNS hierarchy is 
$$u_t=[\p_x-(a\l+u), (Q\l^{j-1})_+]= [\p_x- u, \, Q_{-(j-1)}].$$
By Proposition \ref{cm}, this is an order $j$ evolution partial differential equation. 
For example, the flow generated by $a\l, a\l^2, a\l^3$ for $u=\bpm 0& q\\ r &0\epm$ are
\begin{align*}
&\bca q_t= q_x, \\ r_t= r_x,\eca \quad \bca q_t= \frac{i}{2} (-q_{xx} +2 qrq),\\ r_t=\frac{i}{2} (r_{xx}-2  rqr),\eca \\
& \bca q_t= \frac{1}{4}(-q_{xxx} + 3qrq_x + qr_x q+ 2 q_xrq), \\ r_t=\frac{1}{4}(-r_{xxx} + 3 rqr_x +  rq_x r+ 2r_x q r)\eca.
\end{align*}

\ms

A direct computation implies that  
\beq\label{ih}\tr(vav) =0, \quad \forall\, v\in \calG_1.\eeq 
The following lemma can be proved using \ref{ih} and  induction on $m=i+j$. 

\blem\label{cr}  Let $u=\bpm 0& q\\ r&0\epm$ with $q, r^t\in \C^{1\times n}$, and 
$$\xi_j= \tr(ua^j u^{(j)}), \quad j\geq 0.$$ 
  If $i, j\geq 0$ and $i+j\leq 2n-1$, then
\ben
\item
 $\tr(u^{(i)}u^{(j)})\equiv \tr(u^{(i)}a u^{(j)}) \equiv 0$ {\rm mod} $(\xi_0, \ldots, \xi_{i+j})$,
 \item $\tr(u^{(i)}a^{i+j} u^{(j)})\equiv (-1)^i \xi_{i+j}  \quad {\rm mod\/} \, (\xi_0, \ldots, \xi_{i+j-1})$.
\een
Here we use the notation $A\equiv B$ {\rm mod} $(\xi_0, \ldots, \xi_m)$  if $A-B$ is a polynomial in $\xi_0, \ldots, \xi_m$ and their $t_1$-derivatives.
 \elem
 
Note that $$K:=GL(n+1,\C)_a= GL(n,\C)\times GL(1,\C).$$  Given $f\in L_-(GL(n+1,\C))$,  Write $u_f=\bpm 0& q\\ r&0\epm$, and $k=\bpm k_1&0\\ 0& c\epm\in K$ with $k\in GL(n,\C)$ and $c\in \C$ non-zero. Then 
$$u_{kfk^{-1}}= ku_fk^{-1}=\bpm 0& c^{-1}k_1 q\\ crk_1^{-1}& 0\epm.$$ 
 It is clear that $q^{(i)} r^{(j)}$ is invariant under the $K$-action for all $i, j\geq 0$.  
 
The next lemma proves that these $K$-invariants can be recovered from $\ln\tau_f$. 

\blem\label{ara} Let $u_f=\bpm 0&q\\ r&0\epm$ be the formal inverse scattering solution defined by $f\in L_-(GL(n+1,\C))$ for the vector AKNS hierarchy. Then $q^{(i)} r^{(j)}$ is a  polynomial of $\{y_j=(\ln\tau_f)_{t_1t_k}\n 1\leq k\leq 2n\}$ and their $t_1$ derivatives.
\elem

\begin{proof}
By \eqref{ia}, we have
\beq\label{cq}
(\ln\tau_f)_{t_1t_j}= \li MJ_jM^{-1}, (J_1)_\l\ri_{-1} = \tr(aQ_{-j})=\tr(aT_{-j}),
\eeq
where $T_{-j}$ is given in Proposition \ref{cm}.
We claim that 
\beq\label{cw} y_j:=(\ln\tau_f)_{t_1t_j}= \tr(aT_{-j})\equiv  -\frac{j}{2^j} \xi_{j-1} \quad {\rm mod\,} \, (\xi_0, \ldots, \xi_{j-2}),
\eeq
 where  $\xi_j= \tr(u a^j u^{(j)})$ as in Lemma \ref{cr}. 
To prove this claim, we recall that $va=-av$ for $v\in\calP$ and $$\phi_0=\phi_1=\psi_0= T_0=0.$$  Here $\phi_i$ and $\psi_i$ are given in Proposition \ref{cm}. 
It follows from \eqref{aj2}, \eqref{co1} and \eqref{co2} that we have
\begin{align*}
&-\tr(2aT_{-j})= \sum_{i=0}^{j-1} \tr(P_{-i} P_{-(j-1-i)}+ T_{-i} T_{-(j-1-i)})\\
&= \sum_{i=0}^{j-1} 2^{-(j-1)} \tr\left(a^{j-1} \sum_{i=0}^{j-1} (-1)^i u^{(i)} u^{(j-1-i)}\right) + 2\sum_{i=2}^{j-1}\tr( \phi_iP_{-(j-1-i)}) \\
&\quad  +\sum_{i=2}^{j-3}\tr(\phi_i \phi_{j-1-i})+ \sum_{i=1}^{j-2} \tr(T_{-i} T_{j-1-i})\\
&= A_1+A_2+ A_3+ A_4.\end{align*}
 By Lemma \ref{cr}, we have
$$A_1\equiv \frac{j}{(-2)^{j-1}} \xi_{j-1}\quad {\rm mod\, } (\xi_0, \ldots, \xi_{j-2}).$$ 
To prove \eqref{cw}, it remains to prove that $A_i \equiv 0$ mod $(\xi_0, \ldots, \xi_{j-2})$ for $2\leq i\leq 4$.  This follows from the following two facts:
\ben
\item[(i)] If $v_i=\bpm 0& q_i\\ r_i &0\epm$, then a direct computation implies that
\begin{align*}
\tr(v_1\cdots v_{2k}) &= \tr(r_2q_3)\tr(r_4q_5)\cdots \tr(r_{2k-2} q_{2k-1})\tr(q_1r_{2k})\\ &\quad + \tr(r_1q_2)\tr(r_3q_4)\cdots \tr(r_{2k-1}q_{2k}).\end{align*}
\item[(ii)] By Lemma \ref{cm}, $P_k$, $\phi_k$, $T_k$ are weight $k+1$ polynomials in $u^{(i)}, au^{(i)}$ for $i\leq k$, $i\leq k-2$ and $i\leq k-1$ respectively.  
\een
This proves the claim, i.e., \eqref{cw}. 

A direct computation implies that
\begin{align*} \tr(u^{(i)} u^{(j)}) &= \tr(q^{(i)}r^{(j)}) + \tr(r^{(i)} q^{(j)}) = q^{(i)}r^{(j)} + q^{(j)} r^{(i)}, \\
\tr(u^{(i)} a u^{(j)})&= i(-\tr(q^{(i)}r^{(j)}) + \tr(r^{(i)} q^{(j)}))= i(-q^{(i)}r^{(j)}+ q^{(j)} r^{(i)}.
\end{align*} 
Hence it follows from Lemma \ref{cr}(1) that if $i+j\leq 2n-1$ then
\beq\label{kh}
 \tr(q^{(i)} r^{(j)})\equiv \tr(r^{(i)} q^{(j)})\equiv 0, \quad {\rm mod\/} \quad (\xi_0, \ldots, \xi_{i+j}).
 \eeq  
\end{proof}

We call $h$ a {\it rational differential of $y_1, \ldots, y_k$\/} if $h$ is a rational function of $y_1, \ldots, y_k$ and their $t_1$ derivatives. 

The following is the main result of this section.

\bthm\label{ar} There exist rational differentials  $w_0, \ldots, w_{n-1}, z_0, \ldots, z_{n-1}$ of $y_1, \ldots, y_{2n}$ such that  the formal inverse scattering solution  $u_f=\bpm 0& q\\ r&0\epm$ of the vector AKNS hierarchy defined by $f\in L_-(GL(n+1,\C))$ satisfies
$$\bca
q_{t_1}^{(n)}= w_0(y) q + w_1(y) q^{(1)}+ \cdots + w_{n-1}(y) q^{(n-1)},\\
 r_{t_1}^{(n)}=z_0(y) r + z_1(y) r^{(1)} + \cdots + z_{n-1}(y) r^{(n-1)}.\eca$$
 where $y_i=(\ln\tau_f)_{t_1, t_i}$ for $1\leq i\leq 2n$. \ethm

\begin{proof}
Let $S$ denote the $gl(n, \C)$-valued map whose $i$-th row is $q_{t_1}^{(i-1)}$ for $1\leq i\leq n$, and let $R=(r, r_{t_1}^{(1)}, \cdots, r_{t_1}^{(n-1)})$. We may assume $S$ and $R$ are invertible for generic $u$.   It follows from \eqref{kh} that entries of $C:= SR$ and $b:= q_{t_1}^{(n)}R$ are polynomials of $y_1, \ldots, y_{2n-2}$ and their $t_1$ derivatives, where $y_j= (\ln \tau_f)_{t_1t_j}$.  
So we have  $q_{t_1}^{(n)}=bR^{-1} = bC^{-1} S$.  Hence $q_{t_1}^{n}=W S$, where $W= bC^{-1}$. Similar argument gives the result for $r$.  
\end{proof}

Let $\tau$ be the conjugate involution of $G=GL(n+1,\C)$ defined by $\tau(g)= (\bar g^t)^{-1}$, and $L^\tau_\pm(G)$ the splitting of $L^\tau(G)$ defined by $\tau$ as in Example \ref{bkc}.  Then $\calU= u(n+1)$ is the fixed point set of $\tau$ in $gl(n+1,\C)$ and $\calU= \calK+ \calP$, where 
$$\calK= \left\{\bpm \xi_1&0\\ 0& ir\epm \,\big|\, \xi_1\in u(n), r\in \R\right\},  \quad \calP=\left\{\bpm 0& q\\ -\bar q^t &0\epm\, \big|\, q\in \C^{1\times n}\right\}.$$
The hierarchy constructed from the splitting $L^\tau_\pm(G)$ with $\{a\l^j\n j\geq 1\}$ is the {\it $\C P^n$-NLS hierarchy\/} given in \cite{ForKul83} and the second flow in this hierarchy is the {\it $\C P^n$-NLS equation} (or {\it the coupled $n$ component NLS\/}),
$$q_t=-\frac{i}{2}(q_{xx}+ 2||q||^2 q).$$  Moreover, let $u_f$ and $\tau_f$ denote the formal inverse scattering solution and tau function of the vector AKNS hierarchy defined by  $f\in L^\tau_-(G)$, then $u_f$ is the formal inverse scattering solution and $\tau_f$ is the tau function for the $\C P^n$-hierarchy. 

\bcor There exist rational differentials $w_0, \ldots, w_{n-1}$ of $y_1, \ldots, y_{2n}$ such that the formal inverse scattering solution $u_f=\bpm 0& q\\ -\bar q^t&0\epm$ of the $\C P^n$-hierarchy defined by  $f\in L_-^\tau(GL(n+1,\C))$ satisfies
$$q_{t_1}^{(n)}= w_0(y) q+ w_1(y) q_{t_1}^{(1)} + \cdots + w_{n-1}(y) q_{t_1}^{(n-1)},$$
where $y_i=(\ln\tau_f)_{t_1, t_i}$ for $1\leq i\leq 2n$.  
\ecor

\ms

\brem
Let $U$ be the compact real form of a complex simple Lie group $G$ defined by the conjugate linear involution $\tau$ of $G$, and $a\in \calU$ satisfying $\ad(a)^2=-\id$ on $\calP=\calK^\perp$, where $\calK= \calU_a$. Then $\frac{U}{K}$ is a compact Hermitian symmetric space. Let $\calL_\pm^\tau(\calG)$ denote the splitting of $\calL^\tau(\calG)$ given in Example \ref{bkc}. The second flow in the hierarchy constructed from this splitting and $\{a\l^j\n j\geq 1\}$ is the $\frac{U}{K}$-NLS equation defined in \cite{ForKul83} for maps $u:\R^2\to \calP$.  Hence we call this hierarchy the $\frac{U}{K}$-NLS hierarchy. The calculations given in the proof of Theorem \ref{ar} works for these hierarchies. In particular, the formal inverse scattering solution $u_f$ defined by  $f\in L^\tau_-(G)$ of the $\frac{U}{K}$-NLS hierarchy can be solved from a linear system of ordinary differential equations whose coefficients are rational functions in $\{(\ln\tau_f)_{t_1t_j}\n 1\leq j\leq 2m\}$ and their $t_1$-derivatives, where $m=\dim(\calP)$. 
We also note that the $\frac{U}{K}$-NLS equation is gauge equivalent to the Schr\"odinger flow on $\frac{U}{K}$ (\cite{TerUhl06}).
 \erem

\bs

\section{Tau functions and Virasoro action for the $GL(n,\C)$-hierarchy}\label{pc}

Let $L_\pm$ be the standard splitting of $L=L(GL(n,\C))$, $\calL_\pm$ and $\calL$ the corresponding Lie algebras, and $a=\diag(c_1, \ldots, c_n)$ with distinct $c_1, \ldots, c_n$. The hierarchy constructed from $L_\pm$ and the vacuum sequence 
$$\{a^i\l^j\, \n\, 1\leq i\leq n, j\geq 1\}$$
is the $GL(n,\C)$-hierarchy. The phase space of the flows is $C^\infty(\R, gl(n)_\ast)$, where 
$$gl(n)_\ast=\{(\xi_{ij})\in gl(n,\C)\n \xi_{ii}=0, 1\leq i\leq n\}.$$
Let $\pi_1, \pi_0$ be the projection of $gl(n,\C)$ defined by
$$\pi_1(y)= \sum_{1\leq i\not=j\leq n} y_{ij} e_{ij}, \quad \pi_0(y)=\sum_{i=1}^n y_{ii} e_{ii}$$
for $y=(y_{ij})$. 
In this section, we give a relation between $u_f$ and $\tau_f$ and write down the Virasoro vector fields as partial differential operators for the $GL(n,\C)$-hierarchy.

Let $V(s)= \exp(\sum_{i, j=1}^{n, N} s_{i,j} a^i \l^j)$. Given $f\in L_-$, factor $V(s)f^{-1}= M(s)^{-1}E(s)$. Then 
\beq\label{ig}
u_f(s)= (M(s)a\l M(s)^{-1})_+- a\l
\eeq 
is a solution of the flows in the $GL(n,\C)$-hierarchy generated by $a^i\l^j$ for $1\leq i\leq n$ and $1\leq j\leq N$.  
Next we make a linear change of coordinates for $s_{ij}$'s to make our calculation simpler as follows: 
$$V(s)= \exp(\sum_{i,j=1}^{n, N} s_{i,j} c_k^i e_{kk} \l^j )=\exp(\sum_{k,j=1}^{n, N} t_{k,j} e_{kk} \l^j) = V(t).$$
In other words, 
$$t_{k,j}= \sum_{i=1}^n s_{i,j} c_k^i.$$
Since $\ad(a)$ is a linear isomorphism of $gl(n)_\ast$, there is a unique $v_f\in gl(n)_\ast$ such that
$$u_f(s)= [v_f(t), a].$$
In fact, if $u_f=(u_{ij})$ and $v_f= (v_{ij})$, then 
$$u_{ij}= -(c_i-c_j) v_{ij}.$$
So we have $V(t) f^{-1}= M(t)^{-1}E(t)$. By \eqref{ig}, 
\beq\label{ib}
v_f(t)=\pi_1(m_{-1}(t)),
\eeq  where $m_{-1}(t)$ is the coefficient of $\l^{-1}$ of $M(t,\l)$.

In the following theorem we give a simple relation between $u_f$ and $\ln\tau_f$ written in $t$ variables. 

\bthm\label{ex} Let $u_f=(u_{ij})$ be the formal inverse scattering solution, $v_f=(v_{ij})=[u_f,a]$,  and $\tau_f$ the tau function of the $GL(n,\C)$-hierarchy defined by $f\in L_-(GL(n,\C))$. Then
\beq\label{hk}
 (\ln\tau_f)_{t_{i,1} t_{k,1}} = -v_{ik}v_{ki}, \quad 1\leq i\not=k\leq n. 
\eeq
Or equivalently, 
$$(\ln\tau_f)_{t_{i,1} t_{k,1}}= (c_i-c_k)^2 u_{ik} u_{ki}$$
for $1\leq i\not=k\leq n$, where $a=\diag(c_1, \ldots, c_n)$.
\ethm

\begin{proof} Let $V(t)f^{-1}= M(t)^{-1}E(t)$, 
$$M= \I+ m_{-1}\l^{-1} + m_{-2}\l^{-2} + \cdots,$$  and $v:=v_f(t)= (m_{-1}(t))^\perp$ as above. Set
$$Q_i=Me_{ii}\l M^{-1}= \sum_{j\leq 1} Q_{i,j} \l^j.$$
So $Q_{i,1}= e_{ii}$ and $Q_{i,0}= [v, e_{ii}]$.   

By Theorem \ref{cs}, 
$$(\ln\tau_f)_{t_{i,1} t_{k,1}}= \li Me_{ii}\l M^{-1}, \l e_{kk}\ri_0= \tr(Q_{i, -1} e_{kk}),$$
which is equal to the $kk$-th entry of $Q_{i,-1}$.  
Since $Q_i= Me_{ii} M^{-1}$, we have
$$Q_i = e_{ii}\l + [v, e_{ii}] + Q_{i, -1}\l^{-1} + \cdots.$$ 
To compute $Q_{i, -1}$, we first  note that 
$$M(\I_n -2e_{ii})\l M^{-1}= \l \I_n - 2 Q_i.$$
Since $(\I_n-2e_{ii})^2=\I_n$, we have $(\l\I_n -2 Q_i)^2=\l^2 \I_n$, i.e.,
$$((\I_n-2e_{ii})\l - 2Q_{i,0} - 2Q_{i, -1}\l^{-1} - \cdots )^2=\l^2\I_n.$$
Compare constant term of the above equation to get
\beq\label{gf}
Q_{i, 0}^2 - Q_{i, -1} + (e_{ii} Q_{i, -1} + Q_{i, -1} e_{ii})=0.
\eeq
Note that
\ben
\item[(i)] if $i\not=k$, then the $kk$-th entry of $e_{ii} Q_{i, -1}$ and $Q_{i, -1} e_{ii}$ are zero, 
\item[(ii)] $(Q_{i,0})^2= ([v, e_{ii}])^2 =(\sum_{j\not= i} v_{ji} e_{ji} - v_{ij} e_{ij})^2$. 
\een
 So for $i\not=k$, we get $(Q_{i,-1})_{kk} = -v_{ik}v_{ki}$.
 \end{proof} 

Note that $GL(n,\C)_a$ is the diagonal subgroup.  It follows from Theorem \ref{fb} that if $v=(v_{ij})$ is a solution of the $GL(n,\C)$ hierarchy and $k=\diag(k_1, \ldots, k_n)\in GL(n,\C)$ then $k\cdot v= kvk^{-1}$ is also a solution of the hierarchy. A simple computation implies that
 $$(k\cdot v)_{ij}= \frac{k_i}{k_j} v_{ij}.$$
 Note that $v_{ij}v_{ji}$ is invariant under the action of the diagonal subgroup on the solutions of the $GL(n,\C)$-hierarchy. Hence formula \eqref{hk} agrees with Theorem \ref{fb}.

 \ms
 \beg {\bf [Restrictions of $GL(n,\C)$-hierarchy]}\ 
 
 Let $\sigma$ and $\tau$ be involutions of $G=GL(n,\C)$ defined by 
\beq\label{ic} \sigma(g)=(g^t)^{-1}, \quad \tau(g)= \bar g.\eeq
Then $O(n,\C)$ and $GL(n,\R)$ are the fixed point sets of $\sigma$ and $\tau$ respectively.  Let \
\begin{align*}
&Y=\{(y_{ij})\in gl(n,\C)\n\, y_{ii}=0, 1\leq i\leq n\}, \\
&Y_1=\{(y_{ij})\in Y\n\, y_{ij}= y_{ji}\},\quad
Y_2=\{(y_{ij})\in Y\n\, y_{ij}=y_{ji} \in \R\}.
\end{align*}
The flows in the $GL(n,\C)$-hierarchy are evolution equations on $C^\infty(\R, Y)$ and they leave $C^\infty(\R, Y_1)$ and $C^\infty(\R, Y_2)$ invariant. Moreover, the formal inverse scattering solution $u_f$ for the $GL(n,\C)$-hierarchy defined by $f\in L_-(G)$ lies in $Y_1$ ($Y_2$ resp.) if $f$ is in $ L_-^\sigma(G)$ ($L_-^{\tau,\sigma}(G)$ resp.).  The restrictions of the $GL(n,\C)$ hierarchy to $C^\infty(\R, Y_1)$ and $C^\infty(\R, Y_2)$ are the $(GL(n,\C),\sigma)$-hierarchy and the $\frac{GL(n,\R)}{O(n)}$-hierarchy given in Examples \ref{bkc} and \ref{bla} respectively.  

A formal inverse scattering solution $u_f$ for the $GL(n,\C)$-hierarchy is {\it scaling invariant\/} if $ru_f(r\cdot t)= u_f(t)$ for all $r\not=0$, where $r\cdot t_{i,j}= r^j t_{i,j}$. The scaling invariant solutions of the $(GL(n,\C), \sigma)$-hierarchy are related to semi-simple Frobenius manifolds (cf. \cite{Dub96}).  The $\frac{GL(n,\R)}{O(n)}$-hierarchy arises naturally from the study of flat Egoroff metrics and flat Lagrangian submanifolds of $\C^n$ with flat and non-degenerate normal bundle (cf. \cite{TerWan08}).  
  \eeg
 
 As a consequence of Theorem \ref{ex} we see that $u_f$ of the $(GL(n,\C), \sigma)$ or the $\frac{GL(n,\R)}{O(n)}$ hierarchy is determined by $\ln\tau_f$ up to sign.

\bthm Let $\sigma$ and $\tau$ be the involutions of $GL(n,\C)$ defined by \eqref{ic}.  Given $f\in L_-^\sigma(GL(n,\C))$ ($L_-^{\tau,\sigma}(GL(n,\C))$ resp.), let $v_f=(v_{ij})$ and $\tau_f$ be the formal inverse scattering solution and tau function given by $f$ in the $GL(n, C)$-hierarchy.  Then $v_f$ is symmetric and is a solution of the $(GL(n,\C),\sigma)$ ($\frac{GL(n,\R)}{O(n)}$ resp.) hierarchy. Moreover, 
\beq\label{ky}
(\ln\tau_f)_{t_{i,1}t_{j,1}}= -v_{ij}^2, \quad 1\leq i\not= j\leq n.
\eeq
\ethm

  Note that $O(n,\C)_a= o(n)_a= Z_2^n =\{\e=\diag(\e_1, \ldots, \e_n)\n \e_i=\pm 1\}$. By Theorem \ref{fba}, $Z_2^n$ acts on $v_f=(v_{ij})$ for the $(GL(n,\C), \sigma)$ and $\frac{GL(n,\R)}{O(n)}$ hierarchies by $\e\cdot v= ((\e_i/\e_j)v_{ij})$.  The right hand side of  \eqref{ky} is invariant under this $Z_2^n$ action. So  \eqref{ky} is consistent with Theorem \ref{fba}.

\ms
In the rest of the section, we calculate the formulas for the Virasoro vector fields for the $GL(n,\C)$-hierarchy
induced from the following Virasoro action on $L_-$,
$$\zeta_j(f)= -\left(\l^{j+1} \frac{\p f}{\p \l} f^{-1}\right)_- f, \quad j\geq -1.$$
 
By Theorem \ref{id}, the induced $\calV_+$-action on the reduced frame $M$ and $\ln\tau_f$ are
\begin{align}
& \zeta_\ell(M)M^{-1}= -(\l^{\ell+1} Ef_\l f^{-1}E^{-1})_-, \label{ez3}\\
& \zeta_\ell(\ln\tau_f)= -\li \zeta_{\ell} (M) M^{-1}, E_\l E^{-1}\ri_{-1}. \label{ez2}
\end{align}
First we give the formula for Virasoro vector fields on reduced frames.

\bprop Let $j\geq -1$, and $V(t)=\exp(\sum_{i,j=1}^{n, N} t_{i,j} e_{ii} \l^j)$, and
\beq\label{ez1} \calJ= \l V_\l V^{-1} = \sum_{i, j=1}^{n, m} j e_{ii} \l^j t_{i,j}. 
\eeq
Then the Virasoro vector fields on the reduced frames of the $GL(n,\C)$-hierarchy are
\beq\label{ez} \zeta_j(M)M^{-1}=  -(\l^{j+1} M_\l M^{-1} + \l^j M\calJ M^{-1})_-.\eeq
\eprop

\begin{proof} 
Since $V(t)=\exp(\sum_{i, j=1}^{n,m} e_{ii} \l^j t_{i,j})$, a direct computation implies \eqref{ez1}.

Take the $\l$-derivative of $E=MVf^{-1}$ to get
\beq\label{cj1}
E_\l E^{-1} = M_\l M^{-1} + MV_\l V^{-1} M^{-1}- Ef_\l f^{-1}E^{-1}.
\eeq
This implies that
\beq\label{cja}
\l^{j+1}E_\l E^{-1}= \l^{j+1} M_\l M^{-1} + \l^j M\calJ M^{-1}- \l^{j+1}E f_\l f^{-1} E^{-1}.
\eeq 
Since $j\geq -1$, $E\in \calL_+$, we have $\l^{j+1} E_\l E^{-1}\in \calL_+$ for all $j\geq -1$ and formula \eqref{ez} follows. 
\end{proof}

Let $$\Res_\l(S)= S_{-1}, \quad {\rm if\,}  S=\sum_i S_i \l^i\in \calL(\calG).$$
Since $v_f=\pi_1(\Res_\l(M))$, \eqref{ez} implies the following corollary.

\bcor 
The Virasoro vector fields on $v_f$ for the $GL(n,\C)$-hierarchy is 
$$\zeta_j(v_f)= -\pi_1(\Res_\l(\l^{j+1} M_\l M^{-1}+ \l^j M\calJ M^{-1})).$$
\ecor

Next we calculate the Virasoro vector fields on $\ln\tau_f$ and show that they are given by parietal differential operators. 

\bthm \label{do} Write $\calX= \ln\tau_f$.
The Virasoro vector fields on tau functions of the $GL(n,\C)$-hierarchy given by \eqref{ez2} are partial differential operators. In fact, we have
\begin{align*}
& \zeta_{\ell}\calX = \sum_{i, j=1}^{n, m} j t_{i,j} \calX_{t_{i, j+\ell}} - \frac{1}{2} c_\ell(f), \quad \ell=-1, 0, 1,  \\
 &\zeta_{\ell}\calX = \sum_{i, j=1}^{n,m} j t_{i,j} \calX_{t_{i, j+\ell}} + \sum_{i,j=1}^{n, \ell-1}\left(\calX_{t_{i, j}} \calX_{t_{i, \ell-j}} +\frac{1}{2}  \calX_{t_{i, j} t_{i, \ell-j}}\right) -\frac{1}{2} c_\ell(f), \quad \ell\geq 2.
 \end{align*}
where $c_\ell(f)= \li \l^{\ell+2}(f_\l f^{-1})^2\ri_0$ are constants and $\calX_{t_{i,0}}=0$.
\ethm

\begin{proof}
Use \eqref{ez2} and the fact that $\li \calL_+, \calL_+\ri =0$ to get
\begin{align*}
\zeta_{\ell} \calX&= -\li \zeta_{\ell} (M) M^{-1}, E_\l E^{-1}\ri_{-1}= \li \l^{\ell} E\l f_\l f^{-1} E^{-1}, E_\l E^{-1}\ri_{-1} \\ &= \li \l^\ell E\l f_\l f^{-1} E^{-1}, \l E_\l E^{-1}\ri_0.
\end{align*}
Note that for  for $\xi, \eta\in \calL(gl(n))$ we have
$$\tr(\l^{\ell}\xi\eta)= \frac{1}{2}\left(\l^{\ell}\tr((\xi+\eta)^2)- \tr(\l^{\ell}\xi^2) -\tr(\l^{\ell}\eta^2)\right).$$
So we obtain
\begin{align*}
\zeta_{\ell} \calX
&= \frac{1}{2}\li \l^{\ell} (E\l f_\l f^{-1} E^{-1}+ \l E_\l E^{-1})^2\ri_0-\frac{1}{2} \li \l^{\ell} (E\l f_\l f^{-1} E^{-1})^2\ri_0\\ 
& \quad -\frac{1}{2} \li \l^{\ell}( \l  E_\l E^{-1})^2\ri_0.
\end{align*}
The second term is equal to $-\frac{1}{2} c_{\ell}(f)$, where $c_{\ell}(f)= \li \l^{\ell}( \l^{\ell+2}(f_\l f^{-1})^2\ri_0$ is a constant depending on $f$. The third term is zero because the the lowest order term is $\l^{{\ell} +2}$ and ${\ell}+2\geq 1$ for ${\ell}\geq -1$.  We use \eqref{cja} (with $j=0$) to rewrite the first term and get
\begin{align*}
\zeta_{\ell} \calX &=\frac{1}{2} \li \l^{\ell} (\l M_\l M^{-1} + M\calJ M^{-1})^2\ri_0-\frac{1}{2} c_{\ell}(f)\\
&=\frac{1}{2} \li \l^{\ell}(\l M_\l M^{-1})^2\ri_0 +\frac{1}{2}\li \l^{\ell} \calJ^2 \ri_0 \\
& \qquad + \li\l^{\ell} \l M_\l M^{-1}, M\calJ M^{-1}\ri_0 -\frac{1}{2} c_{\ell}(f) \\
&=\frac{1}{2} (I)+\frac{1}{2} (II)+(III) -\frac{1}{2} c_{\ell}(f),
\end{align*}
where  $\calJ$ is given by \eqref{ez1}. 

Recall that by Theorems \ref{cs} and \ref{csk}, we have
\begin{align}
& (\ln\tau_f)_{t_{i,j}}= \li M e_{ii}\l^j M^{-1}, M_\l M^{-1}\ri_{-1},\label{ag1}\\
& (\ln\tau_f)_{t_{i, j} t_{k,m}}= \li Me_{ii}\l^j M^{-1}, \p_\l(Me_{kk}\l^m M^{-1})\ri_{-1}.\label{ag2}
\end{align}

First we compute (II). 
Since the lowest order term of $\tr(\l^{\ell} \calJ^2)$ is $\l^{{\ell} +2}$ and $\ell\geq -1$, we get
$$(II)=0.$$

We compute (III) as follows:
\begin{align*}
(III)&= \li \l^{\ell} \l M_\l M^{-1}, \sum_{i, j=1}^{n,m} jt_{i, j} Me_{ii} \l^j M^{-1}\ri_0 \\
&=\li \l M_\l M^{-1}, \sum_{i, j=1}^{n,m} jt_{i,j} M e_{ii} \l^{j+\ell} M^{-1}\ri_0, \quad {\rm by\,\, } \eqref{ag1}\\
&= \sum_{i, j=1}^{n,m} j t_{i,j} \calX_{t_{i, j+\ell}}, \quad {\rm if\,\,} \ell \geq 0.
\end{align*}
For $\ell=-1$,  the highest degree term of $\li M_\l M^{-1}, t_{i,1} Me_{ii}\l M^{-1}\ri$ is $\l^{-1}$. So 
$$\li M_\l M^{-1}, t_{i,1} Me_{ii}\l M^{-1}\ri_0=0.$$  Hence by \eqref{ag1} we have
$$(III)= \sum_{i, j=2}^{n,m} j t_{i,j} \calX_{t_{i, j-1}}, \quad {\rm if\,\,} \ell =-1.$$

It remains to compute (I). 
If $\ell\leq 1$, then the degree of  $\l^{\ell+2} (M_\l M^{-1})^2$ is $\ell-2\leq -1$.  Hence 
$$(I)=0, \quad {\rm if\,\, } \ell\leq 1.$$
The computation of (I) for $\ell\geq 2$ is more complicated. Set
\begin{align}
&\xi=(\xi_{ij}):= M^{-1}M_\l. \label{fz1}\\
& \xi_{ij} = \sum_{k\leq -2} \xi_{ij, k}\l^k. \label{fz2}
\end{align}
Note that 
\begin{align*}
(I) &= \li \l^{\ell} (\l M_\l M^{-1})^2\ri_0= \li \l^{\ell+2} (M_\l M^{-1})^2\ri_0\\ 
&= \li \l^{\ell+2} \xi^2\ri_0= \li \l^{\ell+2}\sum_{i, j=1}^n \xi_{ij} \xi_{ji}\ri_0\\
&= \sum_{i=1}^n \li \l^{\ell+2} \xi_{ii}^2\ri_0 + \li \l^{\ell+2} \sum_{i\not=j} \xi_{ij} \xi_{ji} \ri_0\\
&= (I)'+ (I)''.
\end{align*}
By \eqref{ag1}, we get 
\beq\label{gn}
\calX_{t_{i,j}}= \li \l^{j+1} M^{-1}M_\l e_{ii}\ri_0 =\li \l^{j+1}\xi e_{ii}\ri_0= \xi_{ii, -(j+1)}.
\eeq
Thus $$(I)'= \sum_{j=1}^{\ell-1} \calX_{t_{i,j}} \calX_{t_{i, \ell-j}}.$$
To obtain $(I)$, it remains to compute $\sum_{j\not= i} \xi_{ij} \xi_{ji}$.  To do this, we first note that
$$\tr([\xi, e_{ii}]^2)= -2 \sum_{j\not= i} \xi_{ij} \xi_{ji}, \quad {\rm for \,\, } 1\leq i\leq n.$$
Set 
 \begin{align}
 &b_i= \I_n- 2 e_{ii}, \label{fk1} \\
 &P= \l M_\l M^{-1},  \label{fk2}\\
 & Q_i= Me_{ii}\l M^{-1},  \quad B_i= M b_i \l M^{-1}. \label{fm}
 \end{align}
 It is easy to see that $b_i^2=\I_n$,
  \begin{align}
  &[\xi, b_i] = -2 [\xi, e_{ii}], \label{gt}\\
 & B_i^2=\l^2 \I_n, \label{fn0}\\
 &B_i = \l \I_n - 2 Q_i, \label{fn1}\\
 & \l (B_i)_\l = [P, B_i] + B_i, \label{fn2}\\
 & \tr(B_i (B_i)_\l)= n\l. \label{fn3}.
 \end{align}

 Fix $1\leq i\leq n$.  We have
 \begin{align*}
 & (I)''= \li \l^{\ell+2}\sum_{j\not=i} \xi_{ij}\xi_{ji}\ri_0=-\frac{1}{2} \li \l^{\ell+2}[\xi, e_{ii}]^2\ri_0=-\frac{1}{8}\li \l^{\ell+2}[\xi, b_i]^2\ri_0 \\
&= -\frac{1}{8} \li \l^{\ell+2}[M^{-1}M_\l, b_i]^2\ri_0 =  -\frac{1}{8} \li \l^{\ell-2}[\l M_\l M^{-1}, Mb_i \l M^{-1}]^2\ri_0  \\ 
 &=  -\frac{1}{8} \li \l^{\ell-2}[P, B_i]^2\ri_0 =  -\frac{1}{8} (A).
\end{align*}
 We compute $(A)$ next. Recall that $\ell\geq 2$. 
 \begin{align*}
 &(A)= \li \l^{\ell-2} [P, B_i]^2\ri_0, \quad {\rm by\,\,} \eqref{fn2},\\
&= \li \l^{\ell-2} (\l (B_i)_\l - B_i)^2\ri_0\\
&= \li \l^{\ell-2} (\l (B_i)_\l)^2\ri_0+\li \l^{\ell-2} (B_i)^2\ri_0 -2\li \l^{\ell-2} \l (B_i)_\l, B_i\ri_0, \quad {\rm by\,\,} \eqref{fn3}\\
&=\li \l^{\ell} (\p_\l B_i)^2\ri_0 + \li \l^{\ell}\I_n\ri_0- 2\li \l^{\ell} \I\ri_0 \\
& = \li \l^{\ell} (\p_\l B_i)^2\ri_0 =\li \l^{\ell-1} (\p_\l B_i)^2\ri_{-1}\\
&= \li \l^{\ell-1} \p_\l B_i, \p_\l B_i\ri_{-1} = \li \p_\l(\l^{\ell-1} B_i)- (\ell-1)\l^{\ell-2} B_i, \p_\l B_i\ri_{-1} \\
&= \li \p_\l(\l^{\ell-1}B_i), \p_\l B_i\ri_{-1} -(\ell-1)\li \l^{\ell-2} B_i\p_\l B_i\ri_{-1}, \quad {\rm by\,\, }\eqref{fn3}\\
&= \li \p_\l(\l^{\ell-1}B_i), \p_\l B_i\ri_{-1} +(\ell-1)\li \l^{\ell-1}\I_n \ri_{-1} =  \li \p_\l(\l^{\ell-1}B_i), \p_\l B_i\ri_{-1}.
\end{align*}

We need the following formulas for $S, T\in \calL(\calG)$ to compute $(A)$. These formulas can be proved by direct computations.
\begin{align}
&\p_\l (S_+) = (\p_\l S)_+, \quad \p_\l (S_-) = \p_\l (S_-), \label{bk3}\\
&\li \p_\l (ST)\ri_{-1}= \li \p_\l S, T\ri_{-1} +\li S, \p_\l T\ri_{-1}=0, \label{bk4}\\
& \li \p_\l S, T\ri_{-1}= \li \p_\l (S_+), T\ri_{-1}+ \li \p_\l S , T_+\ri_{-1}, \label{bk2}\\
& \li (\p_\l (\l S)_+, T_\l\ri_{-1}= \li \p_\l (S_+)_\l, \p_\l (\l T )\ri_{-1} + 2\li S_+, \p_\l T\ri_{-1}. \label{bk1}
\end{align}

By \eqref{bk2} and \eqref{bk3},
$$(A)=\li \p_\l(\l^{\ell-1}B_i)_+, \p_\l B_i\ri_{-1}+ \li \p_\l(\l^{\ell-1}B_i), \p_\l (B_i)_+\ri_{-1}.$$
The second term is zero because $\p_\l(\l^{\ell-1}B_i)$ has no $\l^{-1}$ term and $\p_\l (B_i)_+= b_i$ is a constant. In fact, we have
$$\li \p_\l X, \p_\l (B_i)_+\ri_{-1}=0$$
for all $X\in \calL$.
 Therefore
 $$(A)= \li \p_\l(\l^{\ell-1}B_i)_+, \p_\l B_i\ri_{-1}.$$
Now we can apply \eqref{bk1} $\ell-1$ times to get 
\begin{align*}
&(A)= \li \p_\l(\l^{\ell-2} B_i)_+, \p_\l (\l B_i)\ri_{-1} + 2\li (\l^{\ell-2} B_i)_+, \p_\l B_i\ri_{-1} \\ &
= \cdots \\
&= \li\p_\l (B_i)_+, \p_\l (\l^{\ell-1}B_i)\ri_{-1} +2\sum_{k=2}^{\ell} \li (\l^{\ell-k} B_i)_+, \p_\l (\l^{k-2} B_i)\ri_{-1}.
\end{align*}
By \eqref{df}, the first term is zero. Hence we get
$$(A)= 2\sum_{k=2}^{\ell} \li (\l^{\ell-k} B_i)_+, \p_\l (\l^{k-2} B_i)\ri_{-1}.$$
Use \eqref{df}, $\li S, T_\l\ri_{-1}= -\li S_\l, T\ri_{-1}$ and \eqref{fn1} to get
\begin{align*}
(A) &= 2\sum_{k=2}^{\ell} \li (\l^{\ell-k} B_i)_+, \p_\l (\l^{k-2} B_i)\ri_{-1}\\
&=-2\sum_{k=2}^{\ell} \li \p_\l(\l^{\ell-k} B_i)_+, \l^{k-2} B_i\ri_{-1}\\
&= -2\sum_{k=2}^\ell \li \p_\l(\l^{\ell-k+1} \I_n- 2\l^{\ell-k} Q_i)_+, \l^{k-1}\I_n - 2\l^{k-2} Q_i\ri_{-1}\\
&= -2\sum_{k=2}^\ell \li (\ell-k+1)\l^{\ell-k}\I_n -2\p_\l (\l^{\ell-k}Q_i)_+, \, \l^{k-1}\I_n - 2\l^{k-2} Q_i\ri_{-1}.
\end{align*}
There are four terms and we will show that three of them are zero: Recall that $\ell\geq 2$. So we have
\ben
\item[(i)] $\li \l^{\ell-1}\I_n\ri_{-1}=0$.  
\item[(ii)] $ \li \l^{\ell-2} Q_i\ri_{-1} =\li \l^{\ell-2} Me_{ii}\l M^{-1}\ri_{-1}= \li \l^{\ell-1} e_{ii}\ri_{-1} =0$.
\item[(iii)] Since $k\geq 2$,  $\li \p_\l (\l^{\ell-k} Q_i)_+, \l^{k-1}\I_n\ri$ is a polynomial in $\l$. So 
$$\li \p_\l (\l^{\ell-k} Q_i)_+, \l^{k-1}\ri_{-1}=0.$$
\een
Therefore we obtain
\begin{align*}
(A) &= -8\sum_{k=2}^{\ell} \li \p_\l (\l^{\ell-k} Q_i)_+, \, \l^{k-2} Q_i\ri_{-1} \\
&=-8\sum_{k=2}^{\ell} \li \p_\l(Me_{ii}\l^{\ell-k+1}M^{-1})_+, \, Me_{ii} \l^{k-1}M^{-1}\ri_{-1}\\
&= -8\sum_{j=1}^{\ell-1} \li \p_\l( Me_{ii} \l^{\ell-j} M^{-1})_+, \, Me_{ii} \l^j M^{-1}\ri_{-1} \\
&= -8 \sum_{j=1}^{\ell -1} \calX_{t_{i, \ell-j}t_{i, j}}.
\end{align*} 
Combine terms to get the formula for $\zeta_\ell \calX$ when $\ell\geq 2$. 
\end{proof}

\brem For $f\in L^\sigma(G)$, a direct computation shows that $\l^{2j-1} f_\l f^{-1}$ is in $\calL^\sigma(G)$. So the vector fields
$$\eta_j(f)=\frac{1}{2}\zeta_{2j} (f)= -\frac{1}{2}(\l^{2j-1}f_\l f^{-1})_-f, \quad j\geq 0,$$
are tangent to $L^\sigma$. Note that $[\eta_j, \eta_k]=(k-j)\eta_{j+k}$ for $j, k\geq 0$.  
Hence Theorem \ref{do} gives the formula for the action of the subalgebra $\calV_0=\{ \eta_j\n j\geq 0\}$ on tau functions $\tau_f$ for the $(GL(n,\C), \sigma)$ and $\frac{GL(n,\R)}{O(n)}$ hierarchies.
\erem


\begin{thebibliography}{99}

\bibitem{AKNS74}
Ablowitz, M.J., Kaup, D.J., Newell, A.C. and Segur, H., \emph{{T}he inverse
scattering transform - Fourier analysis for nonlinear
problems}, Stud. Appl. Math. \textbf{53} (1974), 249--315

\bibitem{Ad79}
Adler, M., \emph{{O}n a trace functional for formal pseudo-differential
operators and the symplectic structure of the Korteweg-de Vries Type
Equations}, Invent. Math., \textbf{50} (1979), 219--248

\bibitem{AraLeu03}
  Aratyn, K., van de Ler, J.,  \emph{An integrable structure based on the WDVV equations}, Theoret. and Math. Phys. 134 (2003), 14Ð26

\bibitem{DriSok84}
 Drinfel'd, V.G., Sokolov, V.V., \emph{{L}ie algebras and equations of Korteweg-de
Vries type}, (Russian) Current problems in mathematics, \textbf{24} (1984),  81--180, Itogi Naukii
Tekhniki, Akad. Nauk SSSR, Vsesoyuz. Inst. Nauchn. i Tekhn. Inform., Moscow

\bibitem{Dub96}
Dubrovin, B.A., \emph{Geometry of 2D topological field
theories}, Lecture Notes in Mathematics, vol. 1620, Springer-Verlag,1996

\bibitem{ForKul83}
Fordy, A.P., Kulish, P.P., \emph{{N}onlinear Schr\"odinger equations and simple
Lie algebra}, Commun. Math. Phys., \textbf{89} (1983), 427-443

\bibitem{Kac85}
Kac, V.G., \emph{{I}nfinite Dimensional Lie Algebras}, Cambridge University
Press (1985)

\bibitem{Kon92}
Kontsevich, M., \emph{{I}ntersection theory on the moduli space of curves and the matrix airy function}, Comm. Math. Phys. \textbf{147} (1992).
 
 \bibitem{PreSeg86}
 Pressley, A., Segal, G. B., \emph{{L}oop Groups}, Oxford Science Publ., Clarendon Press, 
Oxford, (1986) 

\bibitem{ReySem80}
Reyman, A.G., Semenov-Tian-Shansky, \emph{{C}urrent algebras and non-linear partial differential equations}, Sov. Math., Dokl. \textbf{21} (1980), 630-634

\bibitem{Sat84}
Sattinger, D.H., \emph{{H}amiltonian hierarchies on semi-simple Lie
algebras}, Stud. Appl. Math., \textbf{72} (1984), 65--86

\bibitem{Ter08} 
Terng, C.L., \emph{{G}eometries and symmetries of soliton equations and integrable elliptic systems}, Adv. Stud. Pure Math., \textbf{51} (2008) 401--188

\bibitem{Ter14a}
Terng, C.L., \emph{Integrable dispersive curve flows}, in preperation

\bibitem{TerUhl98}
Terng, C.L., Uhlenbeck, K., \emph{{P}oisson actions and
scattering theory for integrable systems}, Surveys in Differential
Geometry: Integrable systems (A supplement to J. Differential
Geometry), \textbf{4} (1998), 315--402

\bibitem{TerUhl00a}
Terng, C.L., Uhlenbeck, K., \emph{{B}\"acklund transformations
and loop group actions}, Comm. Pure. Appl. Math., \textbf{53} (2000), 1--75

 \bibitem{TerUhl06}
 Terng, C.L., Uhlenbeck, K., \emph{{S}chr\"odinger flows on Grassmannians, Integrable systems}, Geometry, and Topology, AMS/IP Stud. Adv. Math. 36 (2006) 235-256

\bibitem{TerUhl11}
Terng, C.L., Uhlenbeck, K., \emph{{T}he $n\times n$ KdV flows}, J. Fixed Point Theory and its Applications, \textbf{10} (2011) 37-61 

\bibitem{TerUhl13b}
Terng, C.L., Uhlenbeck, K., \emph{{T}au functions and Virasoro actions for the $n\times n$ KdV hierarchy}, preprint

\bibitem{TerWan08}
 Terng, C.L., Wang, E., \emph{{T}ransformations of flat Lagrangian immersions and Egoroff nets},  Asian J. Math.  12  (2008),  99--119
 
 \bibitem{vMo94} 
 van Moerbeke, P., \emph{Integrable foundations of string theory. Lectures on integrable systems}, 163Ð267, World Sci. Publ., River Edge, NJ, 1994. 

\bibitem{Wit91}
Witten, E., 
\emph{{T}wo-dimensional gravity and intersection theory on moduli space}, Surveys in differential geometry, \textbf{1} (1990)  243Ð310

\bibitem{Wil91}
Wilson, G., \emph{{T}he $\tau$-functions of the  $\calG$AKNS
equations}, Integrable systems, the Verdier Memorial,  Progress in Math.,
\textbf{115} (1991), 147--162


\end{thebibliography}
\end{document}